\documentclass{article}
\usepackage{amsmath,amsthm}
\usepackage{amssymb}

\newtheorem{prp}{Proposition}[subsection]

\newtheorem{thm}{Theorem}
\newtheorem{cor}{Corollary}[section]

\begin{document}

\newtheorem{lm}[prp]{Lemma}
\newtheorem{rem}[prp]{Remark}
\newtheorem{df}[prp]{Definition}

\numberwithin{equation}{section}


\newcommand{\pendproof}{\endproof}
\newcommand{\goth}{\frak}

\newcommand{ \N }{{\Bbb N}}
\newcommand{ \A }{{\Bbb A}}
\newcommand{ \C }{{\Bbb C}}
\newcommand{ \Dl }{{\Bbb D}}
\newcommand{ \Oc }{{\mathcal{O}}}
\newcommand{ \Q }{{\Bbb Q}}
\newcommand{ \Prj }{{\Bbb P}}
\newcommand{ \U  }{{\mathcal S}}
\newcommand{ \Z }{{\Bbb Z}}
\newcommand{ \CP }{{{\Bbb C}{\Bbb P}^1}}
\newcommand{ \uV }{{\mathbf{\underline {\goth g}}}}
\newcommand{ \wCta }{\widehat{\otimes}{\Bbb C}[[\tau_{\textbf H}]]}
\newcommand{ \bH }{{\textbf H}}
\newcommand{ \V }{{\mathbf {\goth g}}}
\newcommand{ \Vh }{{\mathbf {\goth g}_{\hbar}}}
\newcommand{ \Vht }{{\mathbf {\goth g}_\hbar\,
\widehat\otimes{\Bbb C}[[t_{\textbf H}]]}}
\newcommand{ \G }{{\Gamma}}
\newcommand{ \p }{{\partial}}
\newcommand{ \db }{{\bar {\partial}}}
\newcommand{ \om }{{\omega}}
\newcommand{ \Om }{{\Omega}}
\newcommand{ \dl }{{\Delta}}
\newcommand{ \dlt }{{\delta}}
\newcommand{ \fr }{{,\,\,\,\mbox{for}\,}}
\newcommand{ \Dh }{{\mathcal D}_{\p/\p\hbar}}
\newcommand{ \Df }{{\mbox {Def}}}
\newcommand{ \ga }{{\gamma}}
\newcommand{ \gam }{{\gamma_{\widetilde{\mathcal M}}}}
\newcommand{ \wg }{{\widetilde\gamma}}
\newcommand{ \wgh }{{\widetilde\gamma_\hbar}}
\newcommand{ \gah }{{\gamma_{\hbar}}}
\newcommand{ \Ct }{{{\Bbb C}[[t_{\textbf H}]]}}
\newcommand{ \tng }{{\mathbf {\goth g}}}
\newcommand{\gtg}{{\frak{g}}}
\newcommand{ \Vt }{{\mathbf {\goth g}\,\widehat\otimes{\Bbb C}[[t_{\textbf H}]]}}
\newcommand{ \Ak }{{\Bbb A^{1}_{k}}}
\newcommand{ \AC }{{\Bbb A^1_{\Bbb C}}}
\newcommand{ \Ah }{{\{\Bbb A^1_{k}\setminus 0\}}}
\newcommand{ \Gk }{{\Bbb G_k}}
\newcommand{ \GC }{{\Bbb G_{\Bbb C}}}
\newcommand{ \uH} {{\underline H}}
\newcommand{ \pE} {E}
\newcommand{ \Fh} {{F}_\hbar}
\newcommand{ \Fo} {F_0}
\newcommand{ \Ac} {{\Bbb A_{\Bbb C}^1}}
\newcommand{ \CV} {{\mathcal V}}
\newcommand{ \dE} {{E}\setminus {F_0}}
\newcommand{ \M} {{\mathcal M}}
\newcommand{ \Mh} {{{\mathcal M}_{\hbar}}}
\newcommand{ \MM}  {{\widetilde{\mathcal{M}}}}
\newcommand{ \gta} {{\mathbf {\goth a}}}
\newcommand{ \Ker} {{\mbox{Ker}}}
\newcommand{ \im} {\mbox{Im}}
\newcommand{ \whr} {,\,\,\,\mbox{where}\,\,}

\newcommand{\Remark}{\begin{rem}}
\newcommand{\kendremark}{\end{rem}}





\title{Generalized periods and mirror symmetry in dimensions $n>3$
\thanks
{Based on the author's Ph.D. thesis accepted December 1998}}
\author{Sergey Barannikov\thanks{Alfred P. Sloan Fellow}}
\date{\em Institut des Hautes Etudes Scientifiques,\\
 Bures-sur-Yvette, 91440, France}
\maketitle


\begin{abstract}

The  predictions of the Mirror Symmetry are extended 
in dimensions $n>3$ and are proven for
 projective complete intersections Calabi-Yau varieties. Precisely,
 we prove 
that the total collection of rational Gromov-Witten
invariants of such
variety can be expressed in terms of certain invariants of
a new  generalization of
 variation of Hodge structures attached to the dual variety.

To formulate the general principles of Mirror Symmetry in 
arbitrary dimension 
it is necessary to introduce the ``extended moduli space
of complex structures'' $\mathcal{M}$. 
We show  that the moduli space $\mathcal{M}$ is 
the base of  generalized variation of Hodge structures.
An analog 
$\mathcal{M}\to\oplus_k H^k(X^n,\C)[n-k]$
of the classical period map is described
  and is shown  to be a local isomorphism.
The invariants of the generalized variations 
of Hodge structures are introduced. It is proven 
that their generating function satisfies the system of 
WDVV-equations exactly as in the case of 
Gromov-Witten invariants.
 
 The basic technical tool utilized is the Deformation theory.

\end{abstract}



\tableofcontents


\section{Introduction}\label{ch:int}

This work is devoted
to the description of  the  collection of
all rational
Gromov-Witten invariants of  Calabi-Yau varieties in arbitrary
dimension
 via the invariants of a certain new generalization 
of variations of Hodge structures attached to the mirror dual varieties.

The first discovery in this direction 
was the striking prediction made by 
Candelas, de la Ossa, Green and Parkes \cite{cogp} for the 
numbers of rational curves on quintic threefold in ${\Bbb P}^4$
in terms of the periods of some ``dual'' family of Calabi-Yau
threefolds.

These numbers of rational curves on quintic are the simplest 
examples of 
rational Gromov-Witten invariants.
According to the theory of Gromov-Witten invariants (see \cite{km}) 
the collection of all rational Gromov-Witten invariants of a projective
algebraic
manifold $Y$ is encoded
in the
generating function
\begin{equation}
\text{(Potential)}\,\,\,\,\,\,\,\,\,\,\,\,
\mathcal{F}(t),\,\,\,t\,\, \in\,\,  H^*(Y,\C)
\end{equation}
considered as series over 
 the semigroup ring $\Q[B]$ where $B$ is the
semigroup of effective one-dimensional algebraic cycles modulo
numerical equivalences.
The total space of the cohomology groups
$H^*(Y,\C)$ is considered here as a complex supermanifold.
It is convenient to choose    a graded basis
$\{\Delta_a\}$  in $H^*(Y,\C)$. Let us denote
by $\{t^{a}\}$ the dual set of linear coordinates and choose
some generic representatives $\{\G_a\}$  of the homology classes dual
to $\{\Delta_a\}$ .
Intuitively, the Taylor coefficients $N(a_1,\dots, a_n;\beta)$
in the series expansion 
\begin{equation}
\mathcal{F}(t)=\sum_{n;a_1,\dots,a_n;\beta}
\frac{1}{n!}N(a_1,\dots,a_n;\beta)q^{\beta}
t^{a_1}\dots t^{a_n}
\end{equation}
count the numbers of algebraic maps $f: C\to Y$ 
where $C$ is a rational
curve with $n$ marked  points $\{x_1,\ldots,x_n\}$ such that
$f(x_i)\in \G_{a_i}$ and $f_*([C])=\beta\in B$.
 However a lot of work is needed in order to 
give the precise definition for these numbers (see \cite{bm}). 

The conjectures of \cite{cogp} were partially extended
to higher dimensions in \cite{bvs} (see also \cite{gmp})
where
the formulas describing hypothetically 
a subset of the Gromov-Witten invariants corresponding to the
restriction of the third derivative of the potential
to the subspace of the second cohomology group
\begin{equation}
\p^3\mathcal{F}(t)\,\,\bigl|_{t\in H^2(Y,\C)}
\end{equation}
were proposed.

The next important achievement was made
by  A.~Givental (\cite{g1})
 who has established
 the conjectures from \cite{cogp} and \cite{bvs}
which allow to
express  $\p^3\mathcal{F}(t)|_{H^2(Y,\C)}$ in terms of the
classical periods associated with the mirror dual family.

We construct a generalization of the classical periods map
in order to find the expression for the whole 
generating function $\mathcal{F}(t)$. In other words
the aim of our work  is  to identify
 the total collection of rational Gromov-Witten 
invariants  for the Calabi-Yau varieties of dimension $n>3$ 
with certain invariants
 coming from a  generalization of 
the theory of variations of Hodge structure on the dual varieties.
We prove the coincidence of the two types of invariants for
projective complete intersections Calabi-Yau varieties and its duals.

In section \ref{ch:ext} the ``extended moduli space 
of complex structures'' $\mathcal{M}$ is introduced.
Let $X$ be a complex manifold with trivial canonical sheaf. Let us denote 
by $J_X$ the corresponding complex structure on the underlying 
$C^{\infty}-$manifold $X_{C^\infty}$.
Recall that the tangent space to the classical moduli space of 
complex structures at a smooth point is identified canonically
with $H^1(X,\mathcal{T}_X)$. According
to Kodaira-Spencer theory 
a complex structure on $X_{C^\infty}$ close to $J_X$ is described by an
element $\ga\in \Om^{0,1}(X,\mathcal{T}_X)$ satisfying
Maurer-Cartan equation
\begin{equation}
\db\ga+\frac{1}{2}[\ga,\ga]=0\label{eq:mc}
\end{equation} 
The correspondence: complex structure $J\to \ga_J$
 can be described as follows.
A complex structure on $X$ may be defined as a 
decomposition $T_{\Bbb R}\otimes\C=T\oplus \overline{T}$ 
of the complexified tangent space
into
the sum of complex conjugate subspaces which constitute formally 
integrable distributions (Newlander-Nirenberg theorem).
A deformation of such decomposition  corresponds to a graph 
of a linear map $\overline{T}\to T$, i.e. an element $\ga\in 
 \Om^{0,1}(X,\mathcal{T}_X)$. The equation (\ref{eq:mc}) is
the condition of the formal integrability of $\overline{T}$.
The elements $\ga,\ga'$ describing the 
equivalent complex structures
are related via the action of the group corresponding to the
Lie algebra $\Om^{0,0}(X,\mathcal{T}_X)$.
The extended moduli space of complex structures $\M$ is described 
similarly by the elements 
\footnote{For a graded object $A$ we denote by $A[k]$ the tensor
product of $A$ with the trivial object concentrated in degree $(-k)$}
\begin{equation}
\ga\,\,\,\in\,\,\,\,\oplus_{p,q}\Om^{0,q}(X,
\Lambda^p\mathcal{T}_X)[p-q-1]
\end{equation}
satisfying the eq.~(\ref{eq:mc}). The
technic of the 
deformation theory which allows one to associate
the moduli space $\M$ with the 
differential graded Lie algebra 
$\V=\oplus_{q,p}\Om^{0,q}(X,\Lambda^pT_X)[p-q-1]$
is recalled in \S \ref{ch:ba}.
 A trick from the
 rational homotopy theory 
(see \cite{dgms}) allows one to prove
that the moduli space $\mathcal{M}$
is smooth with the tangent space at the
base point $[X]$ canonically isomorphic to
\begin{equation}
\oplus_{p,q} H^q(X,\Lambda^p
{\mathcal{T}}_X)[p-q]\label{eq:tm} 
\end{equation}
We demonstrate in \S \ref{s:ainfty} using the formality theorem from
\cite{k1} that the supermoduli space
$\mathcal{M}$ parametrizes the
 $A_\infty-$deformations of $D^bCoh(X)$.

The section \S \ref{ch:genper} is devoted to the
descrpition of the generalized period map.
The condition $c_1(T_X)=0$ implies that there exist
nonvanishing holomorphic $n-$form $\Om_{X_t}$, $n=\mbox{dim}_{\C}X$
for every $t\in \M^{classical}$. 
If one fixes a hyperplane $L\subset H^n(X,\C)$ 
transversal to the last component of the Hodge filtration
$F^{\geq n}$ at the base point $[X_0]$
then it allows one
to define the classical period map $\M^{classical}\to  H^{n}(X,\C)$.
It sends a point $(t)$ of the moduli space of complex structures
to the cohomology class of the holomorphic $n-$form $\Om_{X_t}^L$ normalized
so that $[\Om_{X_t}^L]-[\Om_{X_0}]\in L$.
The theorem \ref{th:mainth} and the proposition \ref{pr:piW}
describe a generalization of this map for the extended moduli space
$\M$. Here $W$ is an increasing filtration on the total sum
of cohomology groups $H^*(X,\C)$ complementary to the
Hodge filtration.
It turns out 
that our generalized period map
\begin{equation}
\Pi^W:\M\to \oplus_k 
H^k(X,\C)[n-k],\,\,\,n=\mbox{dim}_{\C} \,X 
\end{equation} is locally an isomorphism.

The  map  $\Pi^W$ arises from  a structure
which may be understood as certain
generalization of the 
variation of Hodge structures (genaralized VHS) having 
the moduli space $\mathcal{M}$ as the base. 
This is explained in sections 
\ref{ch:inv}. 
We introduce also in this section the
invariants of such  
 generalized  VHS.
One of the important properties of these
invariants is the fact that their generating function 
satisfies the system of WDVV-equations
exactly as in the case of Gromov-Witten invariants.

In the section \ref{ch:ms}  we prove that
the rational Gromov-Witten invariants of
projective complete intersections Calabi-Yau manifolds 
coincide with the invariants of the generalized
VHS introduced 
in section \ref{ch:inv} which correspond to their mirror pairs.
One reformulation of this result is the
equality 
$$C_{\alpha\beta}^{\gamma}(\tau)=
\sum_i((\p\Pi)^{-1})^\gamma_i\p_{\alpha}\p_{\beta}\Pi^i$$
where $\p^3_{\alpha\beta\gamma}\mathcal{F}=\sum_{\delta}
g_{\gamma\delta}C_{\alpha\beta}^{\delta}(\tau)$ is the third derivative
of the generating function for rational Gromov-Witten
invariants and $\Pi^i(\tau)$ is the vector of generalized periods
depending on the point of the extended moduli space $\mathcal{M}$
associated with the mirror dual variety.

The importance of the problem of constructing a moduli space 
with the properties 
similar to $\mathcal{M}$
for understanding the 
Mirror Symmetry phenomena was anticipated   
in \cite{w1}. It  was conjectured in \cite{k2}
that such a moduli space is related with 
hypothetical moduli space of $A_\infty-$deformations of $D^bCoh(X)$.
The definition of the moduli space $\M$ appeared
for the first time in \cite{bk}. The relations of our results
with homological mirror symmetry conjecture are discussed in
\S \ref{s:ainfty} and at the end of \S \ref{s:afstr}.

{\sl Acknowledgements.}\,
I am thankful to  my dissertation adviser Maxim Kontsevich
for support and constant interes in my work. I would like also
to acknowledge the stimulating atmosphere of the Institute des Hautes
Etudes Scientifiques in Bures-sur-Yvette where this work was conducted.
My research is supported by the Alfred P. Sloan Fellowship.

All the moduli spaces which we consider are
$\Z-$graded manifolds, in other words
they are supermanifolds with additional 
$\Z-$grading on the structure sheaf compatible with the $\Z_2$-grading.  
To simplify notations we replace $\mbox{deg}\,t^a$ by $\overline{a}$
in superscripts.

\section{Basics of deformation theory}\label{ch:ba}

We recall here the basics of the Deformation theory which
is the main technical tool  used throughout the text.
The material presented here is  well-known
to the specialists  (see for example \cite{schs},\cite{gm},\cite{k1}). 
The 
 Deformation theory was developed in  the
work of a number of mathematicians (P.Deligne, 
V.Drinfeld, B.Feigin, W.Goldman and J.Millson, A.Grothendieck, 
M.Kontsevich,
M.Schle\-ssinger, J.Stasheff \dots). 
Unfortunately many of their results
remained unpublished for a long time. 

We assume that we work over a field~$k$ of characteristic zero.

\subsection{Moduli spaces via differential graded Lie 
algebras}\label{ss:defth}

The principal strategy of Deformation theory may be described
as  follows.
Given some mathematical structure \footnote{for example $A$ can be  an 
associative algebra, complex manifold, vector bundle etc.} $A$
one can associate to $A$ 
the  differential
graded Lie 
algebra\footnote{sometimes it is more convenient to work
with  more general notion of 
$L_\infty-$algebra}
 ${\underline{Der}^*}(A)$
defined canonically up to quasi-isomorphisms.
Recall that the  differential  graded Lie algebra
is a graded vector space equipped with
differential and graded skew-symmetric bracket satisfying a list
of axioms
\begin{gather}
\frak{g}=\oplus_k\frak{g}^k,\,\,\,d:\frak{g}^k\to\frak{g}^{k+1},\,\,
d^2=0,\,\,\,\,\,\,\,\,
[\cdot,\cdot]:\frak{g}^k\otimes\frak{g}^l\to \frak{g}^{k+l}
\notag\\
\,\, d[\ga_1,\ga_2]=[d\ga_1,\ga_2]+(-1)^{\overline{\ga_1}}[\ga_1,d\ga_2],
\,\,\,[\ga_2,\ga_1]=-(-1)^{{\overline{\ga_1}}\,\overline{\ga_2}}
[\ga_1,\ga_2]\label{eq:dgla}\\
[\ga_1[\ga_2,\ga_3]]+(-1)^{\overline{\ga_3}(\overline{\ga_1}+\overline{\ga_2})}
[\ga_3[\ga_1,\ga_2]]+(-1)^{\overline{\ga_1}(\overline{\ga_2}+\overline{\ga_3})}
[\ga_2[\ga_3,\ga_1]]=0\notag
\end{gather} 
The correspondence $A\to{\underline{Der}^*}(A)$ may
be viewed as a kind of ``derived functor''(or rather
``derived correspondence'') 
with respect to the standard correspondence $A\to Der(A)$ which  
associates to $A$ its Lie algebra of infinitesimal
automorphisms.
Then the equivalence classes of
deformations of the structure $A$ 
are described in terms of ${\underline{Der}^*}(A)$.
In the standard approach of the deformation  theory 
 one considers inductevely  
the 
 deformations  up to the given order $1,2,\ldots,N,\ldots$.
In other words, the algebras of functions on the standard parameter spaces
of deformations are the Artin algebras with residue field $k$
(in our context they
will  be  $\Z-$graded  generally).
Recall that such an algebra $\frak{A}$ is  
isomorphic to a direct sum $k\oplus\frak{m}$ where $k$ is a copy of the
base field and $\frak{m}$ is a finite-dimensional commutative
nilpotent 
algebra ($\Z-$graded in general). Even more concretly, any Artin algebra 
with the residue field $k$
is isomorphic to an algebra of
the form
$k[t_i]_{i\in S}/I$, where $I$ is
an ideal $I\supset t^N k[t_i]$ and $S$ is a finite
set of (graded) generators.
Then the deformations $\widetilde{A}_{/\frak{A}}$ of the structure
$A$ over an Artin algebra $\frak{A}$ are   described 
by solutions
to Maurer-Cartan equation 
\begin{equation}
d\ga+\frac{1}{2}[\ga,\ga]=0,\,\,\,\,\ga\in(\underline{Der}^*(A)
\otimes\frak{m})^1
\end{equation}
The equivalent deformations $(\widetilde{A}_{/\frak{A}})_1\simeq
(\widetilde{A}_{/\frak{A}})_2$
correspond to the solutions from the same orbit of the group associated with 
the nilpotent Lie algebra $(\underline{Der}^*(A)\otimes\frak{m})^0$.
This Lie algebra acts on the space $({\frak{g}}
\otimes{\frak{m}})^1$ by
\begin{equation}
\alpha\in({\frak{g}}
\otimes{\frak{m}})^0\to\,\,\,\Dot\ga=d\alpha+[\alpha,\ga]
\end{equation}
It is convinient  to introduce  functor $\Df_{\frak{g}}$
associated with a differential graded  Lie algebra 
$\frak{g}$ 
\begin{equation}
\Df_{\frak{g}}(\frak{A})=\{d\ga+\frac{1}{2}[\ga,\ga]=0|\ga\in({\frak{g}}
\otimes{\frak{m}})^1\}/\G^0({\frak{A}})\label{eq:df}
\end{equation}
which acts from the category
of  Artin algebras with residue field
$k$ to the category of 
sets. We will denote $\Df^{\, 0}_{\frak g}$ the corresponding
functor in the more widely known case when only Artin algebras
concentrated in degree $0$ are involved. Sometimes
we  denote the more general functor on $\Z-$graded Artin
algebras by $\Df^{\,\Z}_{\frak g}$.
The description above of the deformations
of $A$ in terms of $\underline{Der}^*(A)$ may be rephrased now by saying
that 
the functor which associates to $\frak{A}$ the set of
equivalence classes of deformations $\widetilde{A}_{/\frak{A}}$ 
is isomorphic to  $\Df_{\underline{Der}^*(A)}$.

Furthermore, in the cases when the actual
moduli space of deformations of $A$ exists, 
the functor $\Df_{\frak{g}}$ is 
equivalent to the functor $\text{Hom}_{continuous}({\hat\mathcal{O}},\cdot)$
where ${\hat\mathcal{O}}$ is 
the pro-Artin\footnote{=projective limit of Artin algebras}
algebra which is equal to the completion of the algebra of functions
on the actual moduli space of deformations of $A$.

Given a differential graded Lie algebra $\frak{g}$ such that
the functor $\Df_{\frak{g}}$ is equivalent to the functor
represented by some
pro-Artin algebra ${\mathcal{O}}_{\frak{g}}$ one can
define the formal moduli space ${\mathcal{M}}_{\frak{g}}$
associated to $\frak{g}$ by proclaiming ${\mathcal{O}}_{\frak{g}}$
to be ``the algebra of functions on ${\mathcal{M}}_{\frak{g}}$''.

The basic tool to deal with differential graded Lie algebras and
formal moduli spaces associated to them is provided by the
theorem on quasi-isomorphisms which is described
in the next subsection. 

\subsection{Equivalence of Deformation functors}\label{s:BasicTh}

We need to recall first the following 
homotopy generalization of the
notion of the morphism between two differential graded Lie algebras.

A sequence of linear maps
\begin{gather}
 {{F}}_1:{\gtg}_1\to{\gtg}_2\notag\\
{{F}}_2:\Lambda^2(\gtg_1)\to\gtg_2[-1]\label{eq:Fn}\\
{{F}}_3:\Lambda^3(\gtg_1)\to\gtg_2[-2]\notag\\
\cdots\notag
\end{gather}
defines an $L_\infty-$morphism of differential $\Z-$graded Lie algebras 
$\gtg_1$ and 
$\gtg_2$ if 
\begin{multline}\label{eq:Linfty}
d{{F}}_n(\ga_1\wedge\ldots\wedge\ga_n)-\sum_i\pm{
{F}}_n(\ga_1\wedge\ldots\wedge d\ga_i\wedge\ldots\wedge\ga_n)=
\\
=\frac{1}{ 2}
\sum_{k,l\geq 1,\,\,k+l=n} \frac{1}{ k!l!} \sum_{\sigma\in S_n}\pm 
[{{F}}_k(\ga_{\sigma(1)}\wedge\ldots\wedge\ga_{\sigma(k)}),
{{F}}_l(\ga_{\sigma(k+1)}\wedge\ldots\wedge\ga_{\sigma(k+l)})]+\\
+\sum_{i<j}\pm 
{{F}}_{n-1}
([\ga_i,\ga_j]\wedge\ga_1\wedge\ldots\wedge\ga_n)
\end{multline}
In particular, the
first map is a morphism of complexes which respects the Lie brackets  up to 
homotopy defined by the second map, which itself respects the Lie brackets
up to higher homotopies and so on.

An $L_\infty-$map ${F}=\{{F}_n\},{F}:\gtg_1\to\gtg_2$
defines a natural transformations of the functors
${F}_*:\Df_{\gtg_1}\to\Df_{\gtg_2}$. 
If $\ga\in(\gtg_1\otimes\frak{m})^1$ is a solution to Maurer-Cartan
equation then 
\begin{equation}
{F}_*(\ga):=\sum_{n=1}^\infty 
\frac{1}{n!}{F}_n(\ga\wedge\ldots\wedge\ga)\end{equation}
is a solution to Maurer-Cartan equation in $(\gtg_2\otimes\frak{m})^1$

Recall that  an $L_\infty-$morphism 
${F}=\{{F}_n\},{F}:\gtg_1\to\gtg_2$
is called  quasi-isomorphism
if its linear part 
${F}_1$ induces an isomorphism of cohomology of complexes
$(\gtg_1,d_1)$ and $(\gtg_2,d_2)$.
\vskip 2truemm

\noindent{\bf Basic Theorem of Deformation Theory.\,\,}{\it
If ${F}=\{{F}_n\}:\gtg_1\to\gtg_2$ is an  $L_\infty-$morphism  from 
$\gtg_1$ to $\gtg_2$ which is a  quasi-isomorphism then
the natural transformation of deformation functors
${F}_*:\Df_{\gtg_1}\to\Df_{\gtg_2}$ is an isomorphism.}

\subsection{Formal manifolds and odd vector fields}
Here we recall the geometric picture of
 the  theory presented above. In particular
we give an interpretation  of the notions of $L_\infty-$morphism and of
 the moduli space $\mathcal{M}_\gtg$ 
 described by  differential graded Lie algebra $\gtg$.

The set of maps ${F}_n:\Lambda^n\V_1\to\V_2[1-n]$
can be thought of as the 
set of Taylor coefficients of a formal map ${F}:\gtg_1[1]\to\gtg_2[1]$
preserving the zero. Namely,
the algebra of formal power series
on a (super) vector space $\gtg[1]$ can be identified with   
 the dual to the free cocommutative
coalgebra  $\mbox{Symm}(\gtg[1])$ cogenerated by $\gtg[1]$.
Then all geometric objects associated with the $\Z-$graded (formal)
manifold $\gtg[1]$ may be described in terms of this coalgebra.
In particular, a 
map of formal manifolds $\gtg_1[1]\to\gtg_2[1]$
corresponds to a coalgebra
morphism $\mbox{Symm}(\gtg_1[1])\to
\mbox{Symm}(\gtg_2[1])$.
A map of free algebras is defined uniquely by 
its restriction to  the set of generators. Similarly,
a map of coalgebras $\mbox{Symm}(\gtg_1[1])\to
\mbox{Symm}(\gtg_2[1])$ 
is defined uniquely by 
its components
${F}_n:S^n(\gtg_1[1])\to\gtg_2[1],\,\,n\geq 0$.
Recall that in the category of
$\Z-$graded vector spaces $S^n(V[1])=\Lambda^n(V)[n]$. 
The 
set of maps (\ref{eq:Fn}) defines naturally
a (formal power series) map of $\Z-$graded manifolds 
$\gtg_1[1]\to\gtg_2[1]$ preserving the origins,
since the corresponding constant term is zero: $F_0=0$. 
The condition (\ref{eq:Linfty}) can be translated into
geometric terms as follows.

 The structure of the differential graded 
Lie algebra on a graded vector space $\gtg$ is interpreted as 
a degree one vector field $Q_{\gtg}, Q_\gtg^2=0$ on 
the vector space $\gtg[1]$
preserving the origin. Namely, the structure maps
\begin{equation}
d:\gtg\to\gtg[1],\,\,\,\,[\cdot,\cdot]:S^2(\gtg[1])\to \gtg[2]\label{eq:d[]}
\end{equation}
are the components of uniquely defined degree one derivation
 $Q_\gtg$ of the coalgebra $\mbox{Symm}(\gtg)$.
It corresponds to the vector field
\begin{equation}
Q(\ga)=d\ga+\frac{1}{2}[\ga,\ga]
\end{equation}
The relations (\ref{eq:dgla}) which are
satisfied by the structure maps (\ref{eq:d[]}) are exactly
equivalent to the condition
\footnote{The differential graded Lie algebras correspond
to  vector fields whose Taylor expansion contains
only linear or quadratic terms. An arbitrary (formal)
vector field $Q,Q^2=0,Q(0)=0$
of degree one on $\gtg[1]$ is equivalent
by definition to the structure of $L_{\infty}-$algebra on $\gtg$} 
$Q^2_\gtg=0$.

The equations (\ref{eq:Linfty}) are
then reformulated as  the single condition 
\begin{equation}
F_*(Q_{\gtg_1})=Q_{\gtg_2}\label{eq:f(Q)=Q}
\end{equation}

This picture implies that the moduli space described by a
differential graded Lie algebra $\gtg$ can be thought
of as a ``nonlinear cohomology" of the operator $Q_\gtg$.
The precise meaning of this 
may be described as follows. One has a  subscheme 
\begin{equation}
"\mbox{Ker}\,Q":=\{\ga\in\gtg[1]|Q(\ga)=0\}
\end{equation} of 
 zeroes of the
vector field $Q_\gtg$.
The vector fields of the form
$[Q_\gtg,\alpha]$ where $\alpha$ is an arbitrary
constant vector field on $\gtg[1]$ of degree $-1$
span a distribution ``$\mbox{Im}\,Q$''. This distribution
 is tangent 
to the subscheme ``$\mbox{Ker}\,Q$'' since $[Q,[Q,\alpha]]=0$. The
moduli space
$\mathcal{M}_\gtg$  corresponds in geometric terms
to the ``quotent" of the submanifold ``$\mbox{Ker}\,Q$'' defined by the 
zeroes of the vector field $Q$ by the distribution ``$\mbox{Im}\,Q$''.
The functor $\Df_\gtg$ describing the
moduli space $\mathcal{M}_\gtg$ is then identified  with
the natural functor describing the ``quotent" 
``$(\mbox{Ker}\,Q/\mbox{Im}\,Q)$''.


\section{Extended moduli spaces of complex structures}\label{ch:ext}

We recall here following
\cite{bk} the definition
of the extended moduli space of complex structures $\M$.
 We also recall the arguments
 showing  that  it is a smooth moduli space with
the tangent space 
$$
T_{[X]}=\oplus{p,q} H^q(X,\Lambda^p \mathcal{T}_X)[p-q]
$$
\subsection{The sheaf of graded Lie algebras}
Let $X,\,\mbox{dim}_{\C}\,X=n$ be a  smooth 
projective
algebraic manifold such that
 $c_1(T_X)\in\mbox{Pic}(X)$ is zero.
We assume  that $X$ is defined over the complex numbers
although most of our constructions are valid for $X$ defined over an 
arbitrary algebraically closed 
field~$k$ of characteristic zero.

Consider the coherent sheaf of $\Z-$graded Lie algebras 
\begin{equation}\uV=\oplus_k 
\uV^k[k],\,\,\,\,\,\uV^k:=\Lambda^{1-k}{\mathcal{T}}_X\end{equation}
endowed with  the Schouten-Nijenhuys bracket.
This bracket is uniquely defined by the following conditions:
\begin{description}
\item[(1)] For $v_1,v_2\in {\mathcal T}_X$ the commutator $[v_1,v_2]$ is the 
standard bracket on vector fields.
\item[(2)] For $v\in {\mathcal T}_X, f\in {\mathcal O}_X$ the commutator 
$[v,f]\in 
{\mathcal O}_X$ is the Lie derivative $Lie_v f$.
\item[(3)] The wedge 
product and the Lie bracket define the structure of an odd 
Poisson algebra(=Gerstenhaber algebra)  on $\uV[-1]$, in other words
\begin{equation}[v_1,v_2\wedge v_3]=
[v_1,v_2]\wedge v_3+(-1)^{(\overline v_1+1)\overline 
v_2}v_2\wedge[v_1,v_3]\end{equation} 
\end{description}
The  technic of Deformation theory 
associates to the sheaf $\uV$ a (formal) moduli space. 

To describe this moduli space for $X_{/\C}$\footnote{In the case of the  
manifold over arbitrary  field of characteristic zero one should work 
with
 ``simplicial  Lie algebra'' of \v Cech cochains.} 
 let us take the Dolbeault resolution 
of $\uV$ and consider the differential graded Lie algebra
\begin{equation} \V=\oplus_k 
\V^k[k],\,\,\,\,\,\V^k:=\oplus_{k=q-p+1}\Omega^{0,q}(X,\Lambda^pT_X)
\end{equation}
endowed with the differential $\db$  and the extension of 
Schouten-Nijenhuys
bracket by the cup-product of differential forms. 
Informally, the moduli space $\mathcal{M}$
associated with $\V$ may be understood as the moduli space of 
solutions  to Maurer-Cartan equation (\ref{eq:mc}) in $\V$ over
 $\Z-$graded bases modulo gauge 
equivalences.

We show in \S \ref{s:ainfty} that 
 the  moduli space associated to $\V$ by the deformation
theory parametrizes the $A_\infty-$deformations
of $D^bCoh(X)$. 
On the other hand, because of the natural
embedding $\mathcal{T}\subset\Lambda\mathcal{T}$ 
the deformations controlled by the sheaf $\Lambda\mathcal{T}$
generalize
 the deformations of complex structures on $X$. 

\subsection{The (formal) moduli space associated with~
$\Lambda\mathcal{T}$}
 To introduce the  moduli space $\mathcal{M}$
 we use the technic of deformation theory explained
 in \S \ref{ch:ba}.  The moduli space is described by
the deformation functor $\Df_{\V}$
associated to 
$\V$ (see \S \ref{ss:defth}). 
The algebra of functions on the moduli space  $\mathcal{M}$
is by 
definition
 the algebra representing the  
 functor $\Df_{\V}$. First we need to introduce the odd
 ``Laplacian'' operator acting on the space $\V$. The behaviour
 of this operator with respect to various
 algebraic structure on the graded vector space 
$\V$ is described by  the 
 Batalin-Vilkovisky formalism (see for example \cite{schw}).
 
\subsubsection*{Odd Laplacian}
It follows from the condition $c_1(T_X)=0$ that there exists an everywhere 
nonvanishing holomorphic $n-$form $\Omega\in \Gamma(X,\Lambda^n T_X^*)$.
It is defined canonically
up to a multiplication by a constant. Let us fix a choice of 
$\Om$.
It induces isomorphism
of complexes $(\Om^{0,*}(X,\Lambda^pT_X),\db)\simeq 
(\Om^{0,*}(X,\Omega^{n-p}),\db)$;
$\ga\mapsto \ga\vdash\Om$.
One can define then the differential $\dl$  on $\V$ by the formula 
\begin{equation}
(\dl\ga)\vdash\Om=\p(\ga\vdash\Om)\label{eq:delta}
\end{equation}
The Lie bracket on $\V$ satisfies
the following identity (Tian-Todorov lemma)
\begin{equation} 
[\ga_1,\ga_2]=(-1)^{\text{deg}\ga_1+1}(\dl(\ga_1\wedge\ga_2) -(\dl 
\ga_1)\wedge\ga_2 -(-1)^{\text{deg}\ga_1+1}\ga_1\wedge\dl\ga_2)\label{eq:TT}
\end{equation}
where $\mbox{deg}\,\ga=q-p+1$ for $\ga\in\Om^{0,q}(X,\Lambda^p
T_X)$.
In particular $\dl$ is a derivation of the 
differential graded Lie algebra structure.
 
\subsubsection*{Diagram of quasi-isomorphisms}\label{ss:diagram}
Denote by $\bH$ the 
 graded vector space
\begin{equation}\bH=\oplus_k\bH^k,\,\,\,
\mbox{dim}\,\bH^k=\sum_{q-p=k}\mbox{dim}\, 
H^q(X,\Lambda^p T_X)\end{equation} Denote  by $\Ct$ 
the
graded algebra of formal power series on $\bH$.
We recall here the proof from \cite{bk}
of the nonobstructedness  of the $\Z-$graded
moduli space associated to $\V$. Analogous arguments
will be applied later for similar differential graded Lie algebras.
It is convinient to fix some choice of a set 
$\{t^a\}$ of linear coordinates on $\bH$.

\begin{prp}\label{pr:D_t}
 The deformation functor associated to $\V$ is 
 isomorphic 
to the functor represented by the algebra $\Ct$. Equivalently,
there exists a versal solution\footnote{
This means that any other solution over arbitrary Artin
algebra (or projective limit of Artin algebras) is equivalent to a 
solution obtained from
this solution via a base change. It is 
also a "minimal"  solution having these properties.
It is common to use the adjective
 ``versal'' in such situations.}
  to the Maurer-Cartan equation 
\begin{equation}
\db{\ga}(t)+
\frac{1}{2}[{\ga}(t),{\ga}(t)]=0
\end{equation}
in formal power series with values in $\V$ 
\begin{equation}
{\ga}(t)=\sum_a\ga_at^a+
\frac{1}{2!}\sum_{a_1,a_2}\ga_{a_1a_2}t^{a_1}t^{a_2} +\ldots\in 
(\Vt)^1\label{eq:ga(t)}
\end{equation} 
\end{prp}
\begin{rem}The solution of the form (\ref{eq:ga(t)}) is versal iff the 
cohomology classes $[\ga_a]$ form a basis of cohomology of the complex 
$(\V,\db)$.
\kendremark
\proof 
The idea is to use the well-known trick from rational homotopy
theory of K\"ahler manifolds
(see \cite{dgms}) and the theorem from \S \ref{s:BasicTh}.
Notice it follows from the equation
(\ref{eq:TT}) that in the following diagram all arrows are the morphisms
of differential graded Lie algebras 
\begin{equation}
(\V,\db)\leftarrow(\Ker\dl,\db)\rightarrow (\Ker \dl/\im{\dl},d:=0)
\label{eq:ker(dl)}
\end{equation}
The $\p\db-$lemma (see \cite{gh}) implies that both arrows are
in fact quasi-isomorphisms.
\pendproof

\begin{cor}
 One can associate to 
$\V$ the smooth (formal)
\footnote{For our purposes it is
sufficient to work on the level of formal manifolds.
In fact the standard technic of Kuranishi spaces may be used
to show that  all our moduli spaces exist on the level
of analytic manifolds, i.e. all the power series 
representing the versal solutions may be chosen to be convergent}
 moduli space ${\M}$, $\widehat{\mathcal{O}}_{\M}\simeq
\Ct$.
The tangent space to ${\M}$ at the base point $[X]$ is 
canonically 
isomorphic to
the $\Z-$graded vector space $\oplus_{q,p} H^q(X,\Lambda^p T_X)[p-q]$.
\end{cor}

\Remark The differential  Lie algebra $\V$ 
and the associated moduli space were first introduced in \cite{bk} 
where $\V$ was considered with a 
slightly
different grading. The $k-$th graded component of the similar 
differential 
graded Lie algebra $\V^{\tilde *}$ defined  in \cite{bk} is 
$\V^{\tilde k}=\oplus_{q+p-1=k}\Omega^{0,q}(X,\Lambda^pT_X)$. 
The corresponding moduli spaces considered as $\Z_2-$graded (formal) 
manifolds are canonically isomorphic
since 
the two gradings agree mod $\Z_2$. 
\end{rem}

\Remark The classical moduli space 
of complex structures $\M^{classical}$ is  naturally a
subspace of $\M$. The classical deformations 
are parametrized by $t\in\M$ such that $\ga(t)$ is 
equivalent to a solution with values in $\Om^{0,1}(X,T_X)$. 
\end{rem} 

\subsection{Moduli space of $A_\infty-$deformations of $D^bCoh(X)$.}
\label{s:ainfty}
Here we sketch 
demonstration of the fact that the supermoduli space $\M$ 
paramet\-rizes the $A_\infty-$defor\-mations of 
$D^b(Coh(X)$.  The notion of 
$A_\infty-$algebra was introduced by J.Stasheff.
The $A_\infty-$categories in the
context of mirror symmetry were considered for the first
time
  in \cite{k2}. The main result which is used in this subsection
  is the
formality theorem from \cite{k1}.

Recall (see \cite{st}) that
the structure of the $A_\infty-$algebra on a graded vector space
$A$ is  defined by odd degree one derivation $M$ of the 
free coassociative coalgebra generated by $A[1]$
\begin{equation}
\mathcal{C}(A[1]):=A[1]\oplus (A[1])^{\otimes 2}\oplus \ldots\notag
\end{equation}
such that $[M,M]=0$.
Degree one derivations of the free coalgebra 
are in one-to-one correspondence with  collections of
their 
components
$m_k: A^{\otimes k}\to A[2-k]$ for
all $k\in \mathbb{N}$. The condition $[M,M]=0$
is translated into the infinite number
of identities:
\begin{equation}
\sum_{k+l=n+1}\sum_{i=0}^{k-1}
\pm
m_k(x_1,\ldots,x_{i},m_l(x_{i+1},\ldots,x_{i+l}),x_{i+l+1},\ldots,x_n)=0
\label{eq:mm}
\end{equation}
The first identity means that $m_1$ is a differential
on the graded vector space $A$. The next identity
says that $m_1$ is a derivation with respect
to the product defined by $m_2$. The third identity is the
associativity of the product $m_2$ up to the homotopy defined by
$m_3$. The next identity is the compatibility of $m_3$
with the product $m_2$ up to higher homotopies and so on.

There exists a homotopy generalization of the notion of 
module over the associative algebra.
\begin{df}
The structure of  module over an $A_\infty-$algebra
$A$ on graded vector space $E$ is a
degree one connection $M^E$ on the free comodule generated by $E$
\begin{equation}
\mathcal{C}(A[1],E):=E\oplus (E\otimes A[1])\oplus \ldots \oplus
(E\otimes (A[1])^{\otimes k})
\oplus
\ldots
\end{equation}
such that $[M^E,M^E]=0$.
\end{df}
Recall that a  connection on a module $\frak{M}$ over 
differential coalgebra 
$(\mathcal{C},d)$ 
is a linear operator $\nabla:\frak{M}\to\frak{M}$ 
compatible with the coalgebra module structure:
$\Delta^{\frak{M}}\circ\nabla=(\nabla\otimes \,Id\pm 
Id\,\otimes \, d)\circ\Delta^{\frak{M}}$.
Again such  connections acting on free modules 
are in one-to-one correspondence with
collections of linear maps $m_i^E:A^{\otimes i}\otimes E\to E[1-i]$
satisfying an infinite number of identities similar to 
(\ref{eq:mm}). For example, if $A$ is in fact a
differential graded associative algebra and 
$m_i^E=0$ for $i\geq 2$, then $E$ is a  differential graded
module over $A$.

Recall that
an $A_\infty-$category  $\mathcal{X}$ is a collection of objects 
$Ob\, \mathcal{X}$ together with 
graded vector spaces of morphisms $\mbox{Hom}^*(E_1,E_2)$
for any pair  $E_1,E_2\in\mathcal{C}$  and 
"higher compositions"
$$m_k(E_0,\ldots,E_k)
:\mbox{Hom}^*(E_0,E_1)\otimes\ldots\otimes\mbox{Hom}^*(E_{k-1},E_k)
\to \mbox{Hom}^*(E_1,E_k)[2-k]$$
defined for any $k\geq 0$
and $E_0,\ldots,E_k\in Ob(\mathcal{X})$. The maps $m_k(E_0,\ldots E_k)$
 should satisfy
 an infinite number of 
 ``associativity up to homotopy'' constraints 
  similar to (\ref{eq:mm}) which
 can be formulated by saying that $\oplus_{i<j}\mbox{Hom}^*(E_i, 
E_{j})$
 should form an $A_\infty-$algebra. 
 In particular $m_1(E)$ is a differential
 on $\mbox{Hom}^*(E,E)$.
  One assumes usually also  
 that the identity morphism $1_E\in \mbox{Hom}^0(E,E)$ is singled
 out and satisfies $m_2(1_E, f)=f, m_{k}(f_1,\ldots, 1_E,\ldots, f_{k-1})=0$
 for any $k\neq 2$. 

 One can show  that for an $A_\infty-$algebra $A$ the 
 collection of all $A_\infty-$modules
 forms an $A_\infty-$category which we will denote
 by $Mod_\infty (A)$. For example the space 
 $\mbox{Hom}^*(E_1,E_2)$ is defined as the space of all
 comodule morphisms \break
$\mathcal{C}(A[1], E_1)\to \mathcal{C}(A[1], E_2)$.
This is the same as arbitrary collection 
of linear maps $f_k:A^{\otimes k}\otimes E_1\to E_2[-k]$ for
$k\geq 0$.
The differential acting on $\mbox{Hom}^*(E_1,E_2)$ is 
given by $m_1(\phi)= M^{E_2}\circ\phi \pm \phi\circ M^{E_1}$.
 
 The deformations of the $A_\infty-$category $Mod_\infty (A)$
 are described by the functor $\Df^{\,\Z}$ associated with the
 differential graded Lie algebra  which 
consists of derivations
 of free coassociative coalgebra generated 
 by $A[1]$ equipped with counit
 \begin{equation}
 {Der}^*(\mathcal{C}_1(A[1])),\,\,\,
 \mathcal{C}_1(A[1])=k\cdot 1\oplus A[1]\oplus\ldots\oplus
  (A[1])^{\otimes n}\oplus\ldots
 \end{equation}
 The differential is the commutator 
 with $M\in {Der}^1(\mathcal{C}(A[1]))\subset  
{Der}^1(\mathcal{C}_1(A[1]))$.
 This is the differential graded Lie algebra of Hochschild cochains
 $C^*(A,A)[1]$.
Let us suppose that we are given a solution to Maurer-Cartan
equation 
$\Gamma(t)\in (C(A,A)[1]\widehat\otimes \frak{M}_R)^1$
where $R=k\oplus\frak{M}_R$ is a (pro)-Artin algebra
assumed to be $\Z-$graded. Then the 
 objects of the deformed category $(\widetilde{Mod_\infty})_{/R}$
can be described as
pairs 
\begin{gather}
(E,\widetilde{M}^E),\,\, E\in Ob(Mod_\infty),\,\,
\widetilde{M}^E=M^E+\G^E,
\\
\G^E\in (Connec_{\G}^*(\mathcal{C}(A[1],E))\widehat\otimes
\frak{M}_R)^1,
[M^E+\G^E,M^E+\G^E]=0\notag
\end{gather}
 Here $\widetilde{M}^E=M^E+\G^E$ is a connection depending on the
parameters $t\in R$ on the
free comodule $\mathcal{C}(A[1],E)$ which lifts the derivation
$M+\G$ acting
on the coalgebra $\mathcal{C}_1(A[1])$. 
 
 One of the  ways to describe the $A_\infty-$category structure
 on $D^bCoh(X)$ is to identify it with  the category
 which is obtained from the
 differential graded category of finitely generated projective 
modules over the algebra 
$\mathcal{A}_X=(\Om^{0,*}(X),\db)$ by taking
the same collection of objects and  by taking  
$H^0$ of the complex $Hom^* (\mathcal{E}_1^*,\mathcal{E}_2^*)$
as the linear space of morphisms between  two objects $\mathcal{E}_1,
\mathcal{E}_2$.
Then one can use the 
  homotopy technic (see for 
example \cite{mp}) to construct
 the  $A_\infty-$category structure
 on $D^bCoh(X)$. The $A_\infty-$category $D^bCoh(X)$ is
$A_\infty-$equivalent to the category of finitely
generated projective modules over $\mathcal{A}_X$.

Let $\V^{\tilde *}$ denotes our differential graded Lie algebra
of Dolbeault forms with coefficients in polyvector fields
which is considered with slightly different grading
 $(\V^{\tilde *})^k=\oplus_{q+p-1=k}\Om^{0,q}(X,\Lambda^p\mathcal{T})$.
Since the two gradings agree $mod \,2$  the corresponding
 supermoduli spaces 
 are the
 same.
It is possible to prove using the formality theorem from \cite{k1}
that the differential graded Lie algebra $\V^{\tilde *}$ is
quasi-isomorphic to the 
Lie algebra  of  local Hochschild cochains on
differential graded algebra $\mathcal{A}_X$.
Therefore equivalence classes of solutions to Maurer-Cartan
equation in $\V$ parametrize  $A_\infty-$deformations
of $D^bCoh(X)$.
Given a solution to Maurer-Cartan equation $\ga(t)\in
(\V^{\tilde *}\widehat\otimes\frak{M}_R)^1$ the 
quasi-isomorphism referred to above gives similar
solution $\Gamma(t)$ in the Lie algebra of local Hochschild
cochains on $\mathcal{A}_X$.
The objects of deformed category $\widetilde{D^bCoh(X)}_{/R}$
can be identified with pairs 
$(\mathcal{E},\widetilde{\nabla}^{\mathcal{E}})$.
 where
$\mathcal{E}=(\mathcal{E}^*,
\delta^{\mathcal{E}})$ is a bounded complex of locally free sheaves
and $\widetilde{\nabla}^{\mathcal{E}}=(\db+\delta^{\mathcal{E}}+
m_2^{\mathcal{E}})+
\G^{\mathcal{E}}$ is a deformed conection on  
free comodule over $\mathcal{C}(\mathcal{A}_{X}[1])$
generated by $\Om^{0,*}(X,\mathcal{E}^*)$ which
depends on $t\in R$ and
satisfies 
$[\widetilde{\nabla}^E,\widetilde{\nabla}^E]=0$.

\section{Generalized periods.}\label{ch:genper}
To a point of classical moduli space $[X]\in \M^{classical}$
one associates a vector of periods of the holomorphic $n-$form
defined up to multiplication by a constant.
We introduce in this section the  periods
associated to the generalized deformations of complex structure
parametrized by $\M$.  The properties of these periods are
investigated 
 in \S
\ref{ch:inv}.
 
\subsection{Classical periods.}
We recall first the notion of periods associated with 
classical complex structures on $X_{C^\infty}$.  
Let us denote by $J_X$ the complex structure on the underlying
$C^\infty$ manifold $X_{C^\infty}$ corresponding to a point
$[X]\in\M^{classical}$.
The condition 
$c_1(T_X)=0$ (as an element
of $Pic(X)$) implies that there exists holomorphic nowhere 
vanishing $n-$form $\Om_X\in\G(X,\Om^n_X)$
(holomorphic volume element). Such form is defined
canonically up to  multiplication by a constant. 
The deformations of the pair $(J_X,\Om_{X})$ are described
by the deformation functor $\Df^{\,0}_{\underline{Der}^*(J_X,\Om_X)}$
associated with
 the differential graded Lie algebra
 \begin{gather}
 \underline{Der}^*(J_X,\Om_X)=\oplus_q\Om^{0,q}(X,\mathcal{T}_X)[-q]
 \oplus_{q'}
\Om^{0,q'}(X,\mathcal{O}_X)[-q'-1]\\
d:=\db+\dl,\,\,\, [\,\,,\,\,]:=\mbox{Schouten-Nijenhuys 
bracket}
\times\mbox{cup-product}\notag
\end{gather}

 \begin{prp}\label{pr:(ga,f)}
The classes of gauge equivalences of solution to Maurer-Cartan equation 
\begin{equation}
(\rho_J,f)\in\break\Om^1(X,\mathcal{T}_X)
\oplus\Om^0(X,\mathcal{O}_X)
\end{equation} 
 define deformations of the pair $(J_X,\Om_{X})$
consisting of the complex 
structure together with the choice of holomorphic volume form.
The deformations of the complex structures are described by $\rho_J$.
The deformed holomorphic $n-$form is 
\begin{equation}
\mbox{exp}(f+\rho_J)\vdash \Om_{X}\label{eq:wnew}
\end{equation}
\end{prp}
\proof 
The Maurer-Cartan equation for $(\rho_J, f)$ can be written
as 
\begin{gather}
\db\rho_J+\frac{1}{2}[\rho_J,\rho_J]=0\\
\dl\rho_J+\db f +[\rho_J,f]=0
\end{gather}
The first equation means that $\rho_J$ defines
a deformation of the complex structure.
If one chooses a local holomorphic coordinates $(z^1,\dots, z^n)$ on 
$X$ then the  differentials of the deformed set of holomorphic coordinates
are given by 
\begin{equation}
d\tilde z^i=dz^i+\sum_{\bar j}(\rho_J)^i_{\bar j}d{\bar z}^{\bar j}
\end{equation}
where $\rho_J=\sum_{i,{\bar j}}(\rho_J)^i_{\bar j}d{\bar z}^{\bar j}
\frac{\p}{\p z^i}$. One has 
$\Om_{X}=a(z)\prod_i dz^i$ in the local coordinates.
The second equation implies now 
that the $n-$form
\begin{equation}
\mbox{exp}(f+\rho_J)\vdash \Om_{X}=e^f a(z)\prod_i d{\tilde z}^i
\end{equation}
is closed.
One can show similarly that the gauge transformations
preserve the class of equivalenses of the complex structure defined
by $\rho_J$ and the cohomology class of (\ref{eq:wnew})
\pendproof
\Remark
The moduli space of pairs $(J_X,\Om_{X})$ form
naturally a linear bundle $\mathcal{L}\to\mathcal{M}^{classical}$
over 
the moduli space of complex structures.
Analogously one can prove that the formal completion 
at the point $(J_X,\Om_{X})$ of the algebra
of functions  on the total space of the linear 
bundle  $\mathcal{L}$ 
is the algebra representing the
deformation functor $\Df^{\,0}_{\underline{Der}^*(J_X,\Om_X)}$.
\end{rem}

The correspondence $[X]\to \C\cdot\Om_X$ defines the 
(projectivization of the)
periods map $\M^{classical}\to\mathbb{P}H^n(X,\C),\,\,\mbox{dim}_{\C}=n$.
If one chooses a value of $\Om_{X_0}$ at some point 
$[X_0]\in\M^{classical}$ and a hyperplane $L\subset H^n(X,\C)$ 
transversal to the line $\C\cdot\Om_{X_0}$ then one
can define
the periods map  with values in $H^n(X,\C)$: 
\begin{equation}
\Pi^L([X])=\Om^L_{X},
\,\,\,\Om^L_X-\Om_{X_0}\in L\label{eq:norm}
\end{equation}
 Given  a set of elements 
$\{{G}^i_{(n)}\}\in H_n(X,\C)$ which form a locally constant
frame  the map $\Pi^L$ is described by the vector 
of periods $(\int_{G_{(n)}^i}
\Om^L)$.

\subsection{The generalized periods map $\M\to \oplus_k 
H^k(X,\C)[n-k]$}\label{s:thegpm}
We introduce in this section the generalized periods map $\Pi^W:
\M \to \break\oplus_k 
H^k(X,\C)[n-k]$ where $W$ is a filtration on
$\oplus_k H^k(X,\C)[-k]$ complementary to the Hodge
filtration (eq.~\ref{eq:compl}). 
In applications to mirror symmetry $W$ will be 
the limiting weight filtration arising
from degenerating family of Hodge structures.
It is shown in proposition  \ref{pr:piW} that $\Pi^W$ is 
a local isomorphism.

\subsubsection*{Family of differential graded Lie algebras.}
Let us consider the 
family of differential graded Lie algebras \footnote{From here the 
differential graded Lie algebra $\V$ and the 
corresponding moduli space ${\M}$ will be sometimes
referred to also  as 
 $\V_0$ and ${\M}_0$}
\begin{gather}
\Vh=\oplus_k 
\Vh^k,\,\,\,\,\,\Vh^k:=\oplus_{k=q-p+1}\Omega^{0,q}(X,\Lambda^
pT_X)\notag\\
d_\hbar:=\db+\hbar\dl,\,\, [\,\,,\,\,]:=\mbox{(Schouten-Nijenhuys 
bracket)}
\times\mbox{(cup-product)}\label{eq:t_h}
\end{gather}

Here we would like to consider $\Vh$ as a differential graded
Lie algebra over $\C[[\hbar]]$.
Notice that it follows from the $\p\db-$lemma (see \cite{gh}) that the
cohomology of the operator $\db+\hbar\dl$ 
acting on $\oplus_{p,q} \Om^{0,q}(X,\Lambda^p T_X)[p-q-1]\widehat
\otimes\C[[\hbar]]$
is a free $\C[[\hbar]]-$module.
The same diagram (\ref{eq:ker(dl)})
 of quasi-isomorphisms proves the  statement similar
to proposition \ref{pr:D_t} in the case of the family of
differential graded Lie algebras $\Vh$.

\begin{prp}
There exists a solution to Maurer-Cartan equation in $\Vh$ 
depending formally on $\hbar$
\begin{gather}
\db{\wgh}(t)+\hbar\dl{\wgh}(t)+\frac{1}{2}
[{\wgh}(t),{\wgh}(t)]=0\notag\\
{\wgh}(t)=\sum_a(\wgh)_at^a+\frac{1}{2!}
\sum_{a_1,a_2}(\wgh)_{a_1a_2}t^{a_1}t^{a_2} +\ldots,\,\,\,
(\wgh)_{a_1\ldots a_k}\in \V\widehat\otimes\C[[\hbar]]
\end{gather}
such that 
$(\wgh)_{a}$ form a basis of the $\C[[\hbar]]-$module of
the 
cohomology of operator $\db+\hbar\dl$ acting on 
$\V\widehat\otimes\C[[\hbar]]$.
\end{prp}
\pendproof

We will call  solution satisfying this condition
``fibered versal solution to Maurer-Cartan equation in $\Vh$".

If $\wg$ is  such a solution then
for arbitrary family of coordinates change   depending formally
on $\hbar$
\begin{equation}t'(\hbar)=t_{(0)}+\hbar t_{(1)}+\dots \notag
\end{equation}
the formula 
\begin{equation}\wg'(t,\hbar)=\wg(t(\hbar),\hbar) 
\end{equation}
gives again a versal fibered solution to Maurer-Cartan 
equation in $\Vh$.
We will see that an arbitrary filtration on $\oplus_k H^k(X,\C)[-k]$
complementary to the Hodge filtration
defines a way to single out  a canonical choice of such solution.

Introduce the rescalling operator 
\begin{equation}
l_\hbar:\sum\ga^{p,q}\to \sum\hbar^{\frac{p+q-2}{2}}\ga^{p,q}
\,\,\mbox{where}\,\, \ga^{p,q}\in\Om^{0,q}(X,\Lambda^pT_X)
\end{equation}
If $\ga\in (\Vh\otimes \frak{m})^1,\,\,\hbar\neq 0$ 
is a solution to the Maurer-Cartan equation depending
formally on $\hbar$
then $l_\hbar(\ga)$ gives a formal family of
solutions to Maurer-Cartan equation in $\V_{\hbar=1}$.

\subsubsection*{Exponential map.}
\label{ss:exp}

\begin{prp}
For a solution to  Maurer-Cartan equation
$\ga \in (\V_{\hbar=1}\otimes 
\frak{m})^1$  one has  
\begin{equation}
\mbox{exp}(\ga)=1+\ga+\frac{1}{2}\ga\wedge\ga+\ldots
\,\,\in \mbox{Ker}\,\db+\dl
\end{equation}
The $(\db+\dl)-$cohomology class of $e^\ga$ is the same
for gauge equivalent solutions.
\end{prp}
\proof
Let us denote by $m_\ga$ the operator of the contraction
 by $\ga\in (\V\otimes\frak{m})^{odd}$ acting on 
$\Om^{*}\otimes \goth A$. The standard formulas  
\begin{equation}
e^{-A}B e^{A}=B+[B,A]+\frac{1}
{2!}[[B,A],A]+\frac{1}{3!}[[[B,A],A],A]+\ldots +
\end{equation}
\begin{equation}
[[\partial, m_{\ga_1}],m_{\ga_2}]=m_{[\ga_1,\ga_2]}
\end{equation}
imply that
\begin{gather}
e^{-m_\ga}\,\p\,e^{m_\ga}=\p+[\p,m_\ga]+\frac{1}{2} 
m_{[\ga,\ga]}\label{eq:epe}
\\
e^{-m_\ga}\,\db\,e^{m_\ga}
=\db+m_{\db\ga}\label{eq:edbe}
\end{gather}
Applying  the sum of  eq.~(\ref{eq:epe}) and eq.~(\ref{eq:edbe}) to
the 
holomorphic $n-$form $\Om_{X_0}$ and multiplying both sides by 
$e^{\ga}$ one gets 
\begin{equation}
(\dl +\db)e^{\ga}=e^{\ga}\wedge(\dl\ga +\db\ga+
\frac{1}{2}[\ga,\ga])\label{eq:(dl+db)e}
\end{equation}
The independence of the cohomology class on the 
choice of solutions in the same class of 
equivalence follows from (\ref{eq:(dl+db)e})
if one puts $\ga'=\ga+\epsilon \alpha,\,\mbox{deg}\,(\epsilon)=1$.
\pendproof
We will  
denote the cohomology of the complex 
\begin{equation}
C^k=\oplus_{q-p=k} \Om^{0,q}(X,\Lambda^p 
T_X)[p-q],\,\,\,\db+\dl:C^k\to C^{k+1}\label{eq:cmplx}
\end{equation}
by $H$. Of course, $H$ is isomorphic to the 
graded vector space $\bH$, 
but we prefer to reserve the special notation for 
the former vector space since it is defined canonically
and the vector space $\bH$ is only defined up to arbitrary
linear isomorphism.

\subsubsection*{Normalized solution to Maurer-Cartan equation.}
\label{ss:normsol}

Let $F$ be the Hodge filtration which we would like to view 
as the following filtration 
 $F^{\ge n}\subset F^{\ge n-\frac{1}{2}}\subset \dots \subset F^{\ge 0}$ 
  on the  
direct sum of all cohomology groups:\footnote{The reason for such a choice 
of indexes on the subspaces of the filtration will become clear from what
follows}
\begin{gather}
F^{\ge r}:=\oplus_{p-q\geq 2r-n}H^{(p,q)}\label{eq:Fgeqr}\\
H^{(p,q)}:=\{[\phi]\in H^{p+q}(X,\C)[-p-q]|\phi\in\Om^{p,q}\},
\,r\in \frac{1}{2}\Z\notag
\end{gather}
Let 
$$
W_{\le 0}\subset W_{\le \frac{1}{2}}\subset\ldots\subset W_{\le n}
$$
be an increasing filtration which is complementary to the Hodge
filtration in the following sense
\begin{equation}\forall \,\,r
\,\,\,\,\,\,H^*(X,\C)=F^{\ge r}\oplus 
W_{\le r-\frac{1}{2}}\label{eq:compl}
\end{equation}
Here a filtration on $H^*(X,\C)$ is by definition
a filtration in the category of $\Z-$graded vector spaces,
i.e. it is the same as a set of filtrations on every graded component 
of $H^*(X,\C)$. As a side remark let us
notice that the whole construction works analogously 
if $W$ is understood as a filtration in the category of $\Z_2-$graded
vector spaces.

Let us assume  that the choice of the holomorphic
volume form $\Om_{X_0}$ corresponding to the complex structure
$X_0$ is fixed.

\begin{thm}\label{th:mainth}
There exists fibered versal solution $\wgh^W$ to Maurer-Cartan equation in 
$\Vh$ which depends 
formally on $\hbar$
\begin{equation}
\db{\wgh^W}(t)+\hbar\dl{\wgh^W}(t)+\frac{1}{2}
[{\wgh^W}(t),{\wgh^W}(t)]=0\notag
\end{equation}
such that
\footnote{A geometric interpretation of this condition
 will be given in \S
\ref{s:bundle}} 
\begin{gather}
\bigl[\mbox{exp}\,(l_\hbar\wgh^W)-1\bigr]\vdash\Om_{X_0}
\in \oplus_r 
W_{\le r}\hbar^{s(r)}\C[[\hbar^{-1}]]\widehat{\otimes}\C[[t_{\bH}]]
\label{eq:einW}\\
\mbox{where}\,\,\,s(r)=-r+n-1\notag
\end{gather}
\end{thm}
\proof
The condition (\ref{eq:compl}) implies that 
one has the direct sum decomposition 
\begin{equation}
 H^*(X,\C)=\oplus_r F^{\ge r}\cap W_{\le r}\label{eq:hfw}
\end{equation}
Let $\{\Delta_a\}$ be a basis in the vector space $H$ of the cohomology
of  complex (\ref{eq:cmplx})
compatible with the direct sum 
decomposition (\ref{eq:hfw}) so that 
\begin{equation}
\Delta_a\vdash \Om_{X_0}\in 
F^{\ge r_a}\cap W_{\le r_a}
\end{equation}
 The 
cohomology of the complex 
(\ref{eq:cmplx}) are isomorphic to the graded vector space
$\bH$ introduced
in \S \ref{ss:diagram}. Let $\{t^a\}$ be a set of
linear coordinates on $\bH$ which form the basis dual to 
$\{\Delta_a\}$.
 Let 
$$\wgh=\sum_a (\wgh)_a t^a+\frac{1}{2}\sum_{a_1 
a_2}(\wgh)_{a_1 a_2} t^{a_1}t^{a_2} +\ldots$$  denotes a fibered versal 
solution to Maurer-Cartan
equation in  $\Vh$. The linear in $t_{\bH}$ terms $(\wgh)_a$
are in $\mbox{Ker}\,\db+\hbar\dl$. 
One can assume without loss of generality 
that $[l_\hbar(\wgh)_a]=\hbar^{s(r_a)}\Delta_a$.
This can be achieved by a 
substitution 
$$t'(\hbar)=a(\hbar)t,\,a(\hbar)\in GL(\bH)[[\hbar]]$$ which is linear in 
$t_\bH$.
Let us consider now the effect of an arbitrary substitution depending 
formally on $\hbar$ with  linear term equal to identity
$$
(t')^a=t^a_{(0)} +\sum_{k\geq 1} t^a_{(k)}\hbar^k,\,\,\,t^a_{(0)}=t^a
$$
Notice that if one puts 
$$
\bigl[\mbox{exp}\,(l_\hbar\wgh^W)-1\bigr]=\sum_a \Delta^a \hbar^{s(r_a)}
\sum_{k=-\infty}
^{k=+\infty} \Phi^a_{(k)}\hbar^k
$$ 
then
\begin{equation}
\Phi^a_{(k)}=t^a_{(k)}+O(t^2)\fr k\geq 0\label{eq:phiak}
\end{equation}
The condition (\ref{eq:einW}) is equivalent to the system of 
equations
\begin{equation}
\Phi^a_{(k)}=0,\fr k>0\label{eq:phiak=0}
\end{equation}
We see from expression (\ref{eq:phiak}) that this system has unique 
solution in formal power series in $t^a_{(0)}$:
$$t^a_{(k)}=f^a_{(k)}(t^a_{(0)}),\,\,k>0$$
The power series 
\begin{equation}
\wgh^W(t)=\wgh(t+\sum_{k=1}^\infty f_{(k)}(t)\hbar^k,\hbar)
\label{eq:wghW}
\end{equation}
satisfies all the conditions of the theorem. 
\pendproof
\begin{prp}\label{pr:well-def}
The formal power series map 
\begin{gather}
\Phi^W=\bigl[\mbox{exp}\,(l_\hbar\wgh^W)-1\bigr]\vdash\Om_{X_0}
\label{eq:phiW}\\
\Phi^W:\M\to 
\oplus_r 
W_{\le r}\hbar^{s(r)}\C[[\hbar^{-1}]]
\notag
\end{gather}
does not depend on the choice of the solution $\wgh^W$ satisfying
the condition of the theorem \ref{th:mainth} and
depends only on the choice of the
filtration $W$ and the 
value of $\Om_{X_0}$.
\end{prp}
\pendproof

\subsubsection*{The map $\Pi^W$.}
\begin{prp}\label{pr:piW}
The power 
series 
\begin{equation}
\Pi^W=\bigl[\mbox{exp}\,(l_\hbar\wgh^W)|_{\hbar=1}\vdash\Om_{X_0}\bigr]
\in \oplus_k H^k(X,\C)[n-k]\widehat\otimes{\hat 
\Oc}_{\M}\label{eq:exp}
\end{equation}
is well-defined and induces 
isomorphism of germs of (formal)
graded manifolds
\begin{equation}
(\M,[X_0])\to (\oplus_k H^k(X,\C)[n-k],\Om_{X_0})
\label{eq:perm}
\end{equation}
\end{prp}
\proof
If one puts
\begin{equation}
\bigl[\mbox{exp}\,(l_\hbar\wgh^W)\bigr]
=1+\sum_a \Delta^a \hbar^{s(r_a)}(\Phi^W)^a_{(0)}+(\Phi^W)^a_{(-1)}\hbar
^{-1}+\ldots\label{eq:suma}
\end{equation}
then $(\Phi^W)^{a}_{(-i)}
(t)\in\frak{M}_{\C[[t_{\bH}]]}^{i+1},\,
\,(\Phi^W)^a_{(0)}=t^a+O(t^2)$
\pendproof
\Remark
The proposition \ref{pr:well-def} implies that the 
map $\Pi^W:\M\to H^*(X,\C)[n]$ defined 
by the series (\ref{eq:exp})
depends only on the choice of  the 
filtration
$W$ and the value of $\Om_{X_0}$.
\kendremark

\begin{prp}\label{pr:PiMcl}
The restriction of the 
generalized periods map $\Pi^W(t)|_{\M^{classical}}$ coincides
with the classical periods $[\Om_{[X_t]}^L]\in H^n(X,\C)[n]$
normalized according to 
(\ref{eq:norm}) using the hyperplane $L=W_{\le n-1}\cap H^n(X,\C)$.
\end{prp}
\proof
Let $(\rho, f_\rho)$ be a solution to Maurer-Cartan equation
in $\V_{\hbar=1}$ corresponding to the deformations 
$(X_\rho,\Om^{L}_{X_\rho})$ of the pair which consists of complex structure
and holomorphic $n$-form. The $n-$form $\Om^L_{X_\rho}$ is
assumed to be normalized so that $\Om^L_{X_\rho}-\Om_{X_0}\in L$.
Notice that $\ga_\hbar=\rho+\hbar f_\rho$ is  the solution
to the Maurer-Cartan equation in $\Vh$ and
$$
\bigl[\mbox{exp}(l_\hbar\ga_\hbar)-1\bigr]\vdash\Om_{X_0}
=\bigl[\mbox{exp}(\rho+f_\rho)-1\bigr]\vdash\Om_{X_0}\in W_{\le 
n-1}\hbar^0
$$
Therefore
it follows from  (\ref{eq:phiak}) and
 (\ref{eq:phiak=0})  that
one can assume without loss of generality that
\begin{equation}
\wg^W(t,\hbar)|_{t\in \M^{classical}}=\rho(t)+\hbar f_{\rho(t)}
\label{eq:wgw=rho+f}
\end{equation}
Hence 
$\Pi^W(t)=[\Om_{X_t}^L]\fr t\in \M^{classical}$.
\pendproof


\subsubsection*{The map $\mbox{Gr}\,\Phi^W$.}
Let us consider the map $\mbox{Gr}\,\Phi^W:\M\to \mbox{Gr}\,W$ 
defined using the expansion (\ref{eq:suma}) as
\begin{equation}
\mbox{Gr}\,\Phi^W(t)=\sum_a(\Phi^W)^a_{(0)}\Delta_a\vdash\Om_{X_0}
\end{equation}
It follows from the 
equation (\ref{eq:phiak}) that it is a local isomorphism 
as well. In applications to mirror symmetry this map will be used
in order to provide a natural set of coordinates on the moduli space 
$\M$. 
In the sequel it will be convinient to distinguish this set
of coordinates associated with the filtration $W$ and denote
them  by $\{\tau_W\}$.
Let us notice that the elements $\Delta_a\vdash
\Om_{X_0}\in F^{\ge r_a}\cap W_{\leq r_a}$ 
project naturally to a basis in $\mbox{Gr}\, W$.
Let $\{\tau_W^a\}$ denote the corresponding dual basis.
The set of coordinates $\{\tau_W\}$ is characterized by the 
property
\begin{equation}
(\Phi^W)^a_{(0)}=\tau_W^a\label{eq:PhiW0}
\end{equation}

It is easy to calculate the restriction of $\mbox{Gr}\,
\Phi^W$ on $\M^{classical}$. 
The restriction of the solution
$\wgh^W$ to $\M^{classical}$ can be assumed to be
of the form (\ref{eq:wgw=rho+f}).
Hence the map $\mbox{Gr}\,
\Phi^W$ sends the complex structure
$[X_\rho]$ to  the class
\begin{equation}
\mbox{Gr}\,\Phi^W|_{\M^{classical}}:[X_\rho]\to[\Om^L_{X_\rho}-\Om_{X_0}]\in
W_{\leq n-1}/W_{\leq n-\frac{3}{2}}\label{eq:coor}
\end{equation}

\Remark 
Let
$\G_0\in H_n(X,\C)$ be a homology class 
orthogonal to the hyperplane $L=W_{\le n-1}\cap H^n(X,\C)$.
Then the condition (\ref{eq:norm})
may be written as
\begin{equation}\int_{\G_0}\Om_{X_\rho}^L=const \end{equation}
Let $\{\G_1^i\}\subset H_n(X,\C)$ is a set of elements which
defines together with $\G_0$ a basis of the 
subspace of $H_n(X,\C)$ annihilating $W_{\le n-\frac{3}{2}}$.
The  restriction of the map $\mbox{Gr}\,\Phi^W$ on ${\M}^{classical}$
associated with the base point $\rho=0$ and the filtration $W$ 
may be written as the vector
\begin{equation}
\mbox{Gr}\,\Phi^W|_{\M^{classical}}=
(\int_{\G^i_1} \Om_{X_\rho}^L-\int_{\G^i_1}\Om_{X_0})
\end{equation}
It  follows from the condition (\ref{eq:compl}) that
\begin{equation}
\mbox{dim}\,\,(W_{\le n-1}\cap H^n)/(W_{\le n-\frac{3}{2}}\cap H^n)=
\mbox{dim}\,\,H^{n-1,1}
\end{equation}
Therefore the map (\ref{eq:coor}) defines a set of
coordinates on ${\M}^{classical}$.
\kendremark

\section[Invariants of generalized variations of Hodge structures]
{Generalized variations of Hodge structures and its invariants}
\label{ch:inv}

\newcommand{\Gr}{\mbox{Gr}\,W}
\newcommand{\Grr}{\mbox{Gr}^*W}
\newcommand{\wCtW}{\widehat{\otimes}{\Bbb C}[[\tau_{W}]]}
\newcommand{\Gm}{{{\Bbb{G}}_m}}
We have seen in the previous chapter that 
one can define a generalization of the classical periods
map on
the extended moduli space 
${\M}$. 
 We investigate the properties of this 
map in this chapter.
Namely we show that the moduli space ${\M}$
is the base of  generalized  variations of Hodge structures
on $H^*(X,\C)$.
We introduce the invariants of this generalized
VHS in \S \ref{s:invVHS}
This is a collection of polylinear $S_k-$symmetric maps of graded 
vector spaces 
\begin{equation}
A^{(k)}_{[X]}:\,\,\,{\mathcal{H}}^{\otimes k}\,\,\to\,\,\mathcal{H},\,\,\,
\,\,\, \mathcal{H}=\oplus_{p,q}H^q(X,\Lambda^{p}T_X)[p-q]\label{eq:inv}
\end{equation}
associated with  filtration $W$ on $H^*(X,\C)$
complementary to the Hodge filtration.
It is shown  in \S \ref{s:frob}
that  if $W$ is isotropic  with respect to the Poincare pairing
 then
these invariants define a structure of the 
Frobenius manifold on $\M$ (in another terminology
they define a solution to Witten-Dijkgraaf-Verlinde-
Verlinde equation).
We prove in chapter \ref{ch:ms} that
the rational Gromov-Witten invariants 
of the projective complete intersections Calabi-Yau 
varieties  coincide with 
the invariants of generalized VHS 
which are attached to the 
 mirror 
varieties.

\subsection{Nonlinear limits of families of affine spaces}\label{s:nonl}
Let $p: \pE\to \Bbb A^1_{k}$ be a regular morphism  of 
(formal) varieties over 
algebraically closed  field~$k,\,\mbox{char}\, k=0$.  Denote by 
$\Fh$ the fiber lying over the point $\hbar\in \Ak$.
Let the fibers $\Fh,\hbar\neq 0$ are equipped with
affine structure\footnote{Remind that  affine structure 
in the analytic setting is an atlas having
 affine transformations as transition 
functions.} 
 such that the corresponding family of 
affine connections has pole of order one at $\hbar=0$ 
\begin{equation}
\nabla(\hbar):{\mathcal T}_{\pE/\Ak}\to 
{\mathcal 
T}_{\pE/\Ak}\otimes\Omega^1_{\pE/\Ak}[\Fo]\label{eq:G(h)}
\end{equation} 
In other words the connection $\nabla(\hbar)$ is represented as 
\begin{equation}d +\sum_{p=-1}^\infty 
{\G^k_{ij}}^{(p)}\hbar^p, \,\,\, {\G^k_{ij}}^{(-1)}\neq 0
\end{equation}for any choice of (formal)
 trivialization $\pE\sim \Fo\times \Ak$ and 
set
of (formal) coordinates on $\Fo$.

\begin{prp} The term ${\G^k_{ij}}^{(-1)}$ with the highest 
order 
of 
pole defines canonical structure of commutative associative algebra on  
$T\Fo$.
\end{prp}
\proof Note that the 
expression $\G_{ij}^k(\hbar)$ transforms
under a fiberwise change of coordinates $x=f(\tilde x,\hbar), 
f=f_0(\tilde x)+f_1(\tilde x)\hbar+\dots$ 
into\footnote{We assume here for simplicity that the fibers
are purely even varieties.} 
\begin{equation}\sum_l -(J^{-1})^k_l{\frac{\p J^l_j}{\p \tilde x^i}} 
+\sum_{m,n,l}(J^{(-1)})^k_m J^l_i J_j^n\G_{ln}^m
\end{equation}
where $J^k_l(\tilde x,\hbar)=\p f^k/\p \tilde x^l$, and $(J^{(-1)})^k_l$ 
is 
the 
inverse matrix to $J^i_j$.
Therefore  the term ${\G^k_{ij}}^{(-1)}(x)$ having the highest order of 
pole
is a well-defined tensor on $\Fo$. The connection $\nabla(\hbar)$ is 
torsionless: 
$T_{ij}^k(\hbar)=\G^k_{ij}(\hbar)-\G^k_{ji}(\hbar)=0$. Hence,  
\begin{equation}{\G^k_{ij}}^{(-1)}= {\G^k_{ji}}^{(-1)}\end{equation} 
  The connection $\nabla(\hbar)$ is also flat  for all $\hbar\neq 
0$ :
\begin{equation}
d\G(\hbar) +\frac{1}{2}[\G(\hbar),\G(\hbar)]=0
\label{eq:dG=0}
\end{equation} 
  Hence, $[{\G}^{(-1)},{\G}^{(-1)}]=0$ . This can be rewritten as 
\begin{equation}\sum_n 
{\G^k_{in}}^{(-1)}{\G^n_{mj}}^{(-1)}=\sum_n 
{\G^k_{mn}}^{(-1)}{\G^n_{ij}}^{(-1)}\end{equation}
The components of the tensor ${\G^k_{ij}}^{(-1)}$ are the structure 
constants 
of 
the commutative
associative multiplication on (the fibers of) $T\Fo$.\pendproof 

\subsection{Tangent sheaf $\mathcal{T}_\M$.} \label{ss:product}
We describe here the  tangent sheaf of $\M$ and recall
  following \cite{bk} how to construct 
  canonical algebra structure on it.

Consider the space $T_{\ga}$ of first-order deformations (over 
$\Z-$graded bases) of a given versal solution 
to the Maurer-Cartan equation $\ga(t)\in 
(\V\,\widehat\otimes\Ct)^1$. The $\Ct-$module 
$T_{\ga}$   is identified naturally with the cohomology of the complex 
$(\V\,\widehat\otimes\Ct[1],
 \db+[\ga(t),\,\cdot\,])$. It is a  $\Z-$graded module over 
the algebra of functions on the moduli space. In geometrical language it 
corresponds to the tangent sheaf of ${\M}$.

The 
operator 
$\db$ as well as the operators
$\mbox {ad}_\ga,\,\ga\in\V$ are differentiations with respect to the 
natural
algebra structure on $\V[1]$ defined by the wedge product.
Hence the $\Ct-$linear extension of the wedge product defines the 
natural $\Oc_{\M}-$linear algebra structure on $T_{\ga}$.
This algebra structure is functorial with respect to the
isomorphisms $T_{\ga_1}\sim T_{\ga_2}$ induced by the  gauge 
equivalences $\ga_1\sim\ga_2$

The 3-tensor  of structure constants of  the multiplication on 
$T_{\ga}$ can be written explicitly as follows. 
Recall that we have fixed for
convinience a choice $\{t^\alpha\}$ of linear coordinates on 
$\textbf H$. 
The 
(uni)versality property of $\ga(t)$ implies that
the classes of partial derivatives $[\p_\alpha\ga]$ generate freely the
$\Ct-$module $T_{\ga}$.
\begin{prp} The equation
\begin{equation}\p_\alpha\ga(t)\wedge\p_\beta\ga(t) = 
\sum_\delta A_{\alpha\beta}^\delta(t)\,\p_\delta
\ga(t) 
\,\,\,\,\mbox{mod Im}\,\,\db_{ \ga(t)}
\label{eq:product}
\end{equation} uniquely determines 
the  formal
power series 3-tensor $A_{\alpha\beta}^\delta(t)\in \Ct$. The components 
of the tensor 
 are 
the  
structure constants of commutative associative \break $\hat\Oc_{\M}-$algebra 
structure.
\end{prp}
 We will denote the product of two elements $v_1,v_2\in \mathcal{T}_{\M}$
 by $v_1\circ v_2$.

Recall that the third derivative of the generating function
for Gromov-Witten invariants of some projective
variety $Y$ defines the commutative  
algebra structure on the tangent sheaf to $H^*(Y,\C)$ considered
as a supermanifold.  We have found 
 similar structure on $\M$. Although, such a structure
alone is relatively noninteresting (for example in many
cases one can choose locally some coordinates $\{u_i\}$ on the underlying 
supermanifold $H^*(Y,\C)$ such that the multiplication is given simply 
by $\frac{\p}{\p u^i}\circ\frac{\p}{\p u^j}=\delta_{ij}\frac{\p}{\p u^i}$),
it indicates the existence of canonical nonlinear family of affine spaces 
 having the former supermanifold as the central fiber.

\subsection{Family of moduli spaces}
Our aim now is to demonstarte that the extended moduli space of 
Calabi-Yau manifolds introduced in chapter \ref{ch:ext}
 is naturally the limiting fiber
of the family of affine spaces with the property (\ref{eq:G(h)}).
We describe first the canonical family of the 
moduli spaces $\Mh$.

\subsubsection*{Extended moduli spaces of complex manifolds with holomorphic
volume element}
In this subsection the moduli space $\M$ 
is embedded into a natural $1-$parameter family of 
 moduli spaces
$\Mh,\,\M_0=\M$.

Let us 
consider a family of sheaves of differential graded Lie 
algebras $(\uV,\hbar\dl)$.
The deformation theory associates to $(\uV,\hbar\dl)$
certain moduli space $\Mh$.
It is convinient to describe this moduli space
with the help of Dolbeult resolution of $(\uV,\hbar\dl)$.
This is  the family of differential
graded  Lie algebras $\Vh$ introduced in \S \ref{s:thegpm}.   
Same diagram (\ref{eq:ker(dl)})
 of quasiisomorphisms proves that  
 the deformation functor associated with the 
differential graded Lie algebra
$\Vh$ for arbitrary given
$\hbar$ is equivalent to the functor represented by the  algebra 
$\Ct$. This can be reformulated as usual as the existence 
of a versal solution to the Maurer-Cartan equation
$\gah(t)\in (\Vh\widehat\otimes\Ct)$
   whose linear term
  gives a basis of  the graded vector space equal to the
  cohomology of the complex $(\Vh,\db+\hbar\dl)$.

According to the proposition \ref{pr:(ga,f)}
the subspace of the moduli space $\M_{\hbar=1}$
corresponding to solutions with values in $\Om^{0,1}(X,T_X)\oplus
\Om^{0,0}(X,\Oc_X)\subset
\V_{\hbar=1}$
is the moduli space of deformations of pair which
consists of complex structure and  
holomorphic $n-$form.
\subsubsection*{The total space of the family $\Mh$}
The family $\Mh, \hbar\in \AC$ of moduli spaces form the 
bundle 
 $p:\MM\to\Ac$. The total space of the bundle can be described 
similarly as the moduli space associated with the differential graded
Lie algebra 
$\widetilde{\V}$. This is the algebra $\V$ with one added element
$D$ of degree $1$ such that $[D,\ga]:=\dl\ga$. Again, the arguments 
involving the use of the diagram (\ref{eq:ker(dl)})
prove that $\widetilde{\V}$ is quasi-isomorphic to an an abelian Lie 
algebra$\widetilde{\textbf H}:={\textbf H}\oplus {\textbf D}$, where 
$\textbf{D}$ is a one-dimensional vector space in degree $0$. Denote by 
$\hbar$  thecorresponding linear coordinate along $\textbf{D}$, 
$\mbox{deg}\,\hbar=0$\begin{prp}\label{pr:D_wt}
 The deformation functor associated with $\widetilde{\V}$
is isomorphic to the functor represented by the algebra 
$\C [[t_{\widetilde
{\bH}}]]$. Equivalently, there exists a versal solution to 
Maurer-Cartan equation in $(\widetilde{\V}\,\widehat\otimes 
t_{\widetilde{\bH}}\C[[t_{\widetilde{\bH}}]])^1$. 
This solution can be taken to be of the form
$\ga_{\MM}=\wgh+\hbar D$ where 
\begin{gather}
 \db{\wgh}(t)+\hbar\dl{\wgh}(t)+\frac{1}{2}
[{\wgh}(t),{\wgh}(t)]=0 \notag\\
\wgh=\sum_a(\wgh)_at^a+\frac{1}{2!}
\sum_{a_1,a_2}(\wgh)_{a_1a_2}t^{a_1}t^{a_2}+\ldots
\label{eq:wgh(t)}
\\ 
\wgh\in 
(\V\,\widehat{\otimes}\C[[t_{\widetilde{\textbf{H}}}]])^1,
\,\,\,(\wgh)_{a_1\ldots a_k}\in \V\widehat\otimes\C[[\hbar]]\notag
\end{gather}
\end{prp}
\proof Define the differential graded Lie subalgebra $\widetilde{ 
\Ker}:=\Ker[D,\,\,]\subset \widetilde{\V}$. It is a 
direct sum of two differential 
graded
Lie subalgebras $\widetilde{\Ker}=D\oplus\Ker\dl,\Ker\dl\subset\V$.
The $\p\db-$lemma implies again that the natural embedding 
$\widetilde{\Ker}\subset \widetilde{\V}$ and the projection 
$\widetilde{\Ker}\to 
\widetilde {\bH}$ are quasi-isomorphisms. 
The versal  solution can be 
taken to be $\ga_0(t)+\hbar D$ where $\ga_0(t)=\sum_n \frac{1}{n!}
f_n (t)$ 
and $f_n:S^n\bH 
\to\Ker\dl$ are the components of a quasi-isomorphism 
$\bH\to \Ker \dl$ 
homotopy inverse to the natural projection $\Ker\dl\to\bH$.
\pendproof 

\begin{rem} A solution of the form (\ref{eq:wgh(t)}) is versal iff
 the cohomology classes $[(\wgh)_a]$ freely 
generate the $\C[[\hbar]]-$module of cohomology of the complex 
$(\V\widehat\otimes\C[[\hbar]],\break\db+\hbar\dl)$
\end{rem}
We see that the  
``fibered versal solutions'' introduced  in \S \ref{s:thegpm}.
describe the total space of the bundle $\MM\to\AC$. 
The choice of class of gauge equivalence of such solutions
corresponds to the fiberwise choice of coordinates
$\MM\to\bH\times\AC$.

\subsection{Affine structure on $\Mh,\hbar\neq 0$}
\label{s:afstr}
  In this subsection
we show that the moduli spaces $\Mh$ for $\hbar\neq 0$
have canonical affine structure. It is induced by the exponential
map from \S \ref{ss:exp}.

Let 
\begin{equation}\gta=\oplus_k 
\gta^k,\,\,\,\gta^k:=\oplus_{k=q-p+1}\Omega^{0,q}
(X,\Lambda^pT_X)
\end{equation}  
denotes the 
differential $\Z-$graded abelian (i.e. with 
zero bracket) Lie algebra   equipped with the differential 
$\db+\dl$. Since 
$\gta$ 
is abelian, the deformation functor associated to 
$\gta$ is representable and the  moduli space associated with 
$\gta$
 is canonically
 isomorphic to (formal neighborhood of zero in)
the cohomology of the complex $(\gta,\db+\dl)[1]$. Therefore 
it has canonical 
structure 
of 
affine
space.

\begin{prp}The differential graded Lie algebra $\Vh$ with 
$\hbar\neq 0$ is 
quasi-isomorphic to the abelian Lie algebra $\gta$.
\end{prp}
\proof 
The differential graded Lie algebras $\Vh$ are all quasi-isomorphic for 
$\hbar\neq 0$. The quasi-isomorphism is given by
the map \begin{equation}i_\hbar:\Vh\to\V_{\hbar=1},\,\,\,\, 
i_\hbar:\phi^{p,q}\to 
 \hbar^{p-1} \phi^{p,q},\,\,\, 
\mbox{for}\,\,\,\phi^{p,q}\in \Omega^{0,q}(X,\Lambda^pT_X)\end{equation}
Therefore it is enough to consider the differential graded Lie algebra
$\V_{\hbar=1}$.
It turns out that one can continue the identity map 
$\phi_{(1)}:\V_{\hbar=1}\to \gta$ of complexes of vector spaces
to an $L_\infty-$morphism of differential graded Lie algebras.

Define the higher components of the $L_\infty-$map 
 $\phi=\{\phi_{(n)}\}:\V_{\hbar=1}\to\gta$ as follows 
\begin{gather}
\phi_{(2)}:\ga_1\wedge\ga_2 \to \ga_1\cdot\ga_2 
\\
\phi_{(3)}:\ga_1\wedge\ga_2\wedge\ga_3 \to \ga_1\cdot\ga_2\cdot\ga_3\\ 
\dots 
\end{gather}
where the symbol "$\cdot$" denotes the  natural wedge 
product on $\gta[1]$.  
Let $\ga_1,\dots,\ga_n\in\V_{\hbar=1}$. Consider the Artin algebra $\goth 
A:=\otimes_{i=1}^n\{\C[\epsilon_i]/\epsilon_i^2=0\}$ where
$\mbox{deg}\epsilon_i+\mbox{deg}\ga_i-1=0$.
Let us look at the equation (\ref{eq:(dl+db)e}) for
 $\ga=\sum_i\ga_i\epsilon_i$. Considering the coefficients 
in front of the 
term $\Pi_i \epsilon_i$ in the equation (\ref{eq:(dl+db)e})
 we see that $\{\phi_{(n)}\}$ has 
the required property (\ref{eq:Linfty})\pendproof 
\Remark  The equation (\ref{eq:(dl+db)e})
 can be directly interpreted 
as the required compatibility of the map $\{\phi_n\}$ of formal 
manifolds with the action of vector fields $Q_{\V_{\hbar=1}}$ and 
$Q_{\gta}$ (see equation (\ref{eq:f(Q)=Q})).
The left hand side is the value of the vector 
field $Q_{\gta}$ at the point
$\phi(\ga)$ and the right hand side is the vector 
field $\phi_*(Q_{\V_{\hbar=1}})$.
\kendremark

\begin{cor} 
The moduli spaces $\Mh$, $\hbar\neq 0$
have canonical affine structure. A family of sets of affine 
coordinates is given by  (compare with formula (\ref{eq:wnew}))
\begin{equation}\int_{G_i} 
 e^{l_\hbar\wgh(t)}\vdash\Om_{X_0}\end{equation}
where $\{G_i\}$ is a basis in $H_*(X,\C)$ and $\Om_{X_0}$ is 
the holomorphic nowhere vanishing $n-$form on~$X_0$.
\end{cor} 

\noindent{\bf Completed tensor products.}
Let ${\mathcal T}_1,{\mathcal T}_2$ are finetely generated 
modules over the topological(=projective limit of Artin) algebras
${\mathcal O}_1,{\mathcal O}_2$. We  denote as usual by 
${\mathcal T}_1\widehat\otimes{\mathcal T}_2$ the completion of 
the tensor product
with respect to the filtration by
 ${\mathcal O}_1\otimes{\mathcal O}_2-$submodules
${\mathcal F}_i=\oplus_{i+j=k} {\goth M}_{{\mathcal O}_1}^i\otimes 
{\goth M}_{{\mathcal O}_2}^j$. 
${\mathcal T}_1\widehat\otimes{\mathcal T}_2$ is 
naturally a ${\mathcal O}_1\widehat\otimes{\mathcal O}_2-$module.
Let $\mathcal{N}_\hbar$ be
 a finitely generated $\C(\hbar)-$module (or $\C((\hbar))-$module)
  and ${\mathcal T}$ be a finetely generated ${\mathcal O}-$module 
 where ${\mathcal O}$ is a
 pro-Artin algebra  with the maximal ideal
${\goth M}_{{\mathcal O}}$. 
The localization 
at $\hbar=0$ of the 
$\C[[\hbar]]\widehat\otimes\mathcal{O}-$module 
$\mathcal{N}_\hbar\widehat\otimes\mathcal{T}$  consists of the elements of 
the form$$\sum_{k=1}^{k=N} \hbar^{-k} \Delta_k, \,\,\,\Delta_k\in
\mathcal{N}_\hbar\widehat\otimes
\mathcal{T}$$
We denote
via  $\mathcal{N}_\hbar\widehat\otimes_{(0)}{\mathcal O}$ the 
minimal completion of the localization at $\hbar=0$ of  
$\mathcal{N}_\hbar\widehat\otimes
\mathcal{T}$ which contains all the elements of the form
$$\sum_{k=1}^{k=+\infty} \hbar^{-k} f_k(t)\,\Delta,\,\,f_k(t)\in 
\frak{M}^k_{\mathcal{O}},\Delta\in 
\mathcal{N}_\hbar\widehat\otimes\mathcal{T}$$ 
In particular $\mathcal{N}_\hbar\widehat\otimes_{(0)}\mathcal{T}$
is a $\C[[\hbar]]\widehat\otimes\mathcal{O}-$module.

\Remark 
The affine structure on $\Mh$ defines the canonical affine map $\Mh\to 
T_{\pi_\hbar}$. It is distinguished from  other
affine maps by the condition that it induces an identity map on 
$T_{\pi_\hbar}$.
In terms of a versal fibered solution $\wgh$ 
it can be written as
\begin{equation}\Phi^H(t,\hbar)=\bigl[\hbar\mbox{exp}\,(\frac{1}{\hbar}
\wgh)-1\bigr]\in 
H_{\hbar}\widehat{\otimes}_{(0)}\C[[t_{\bH}]] \label{eq:phi(t,h)}
\end{equation}where $H_{\hbar}$ 
is the graded $\C((\hbar))-$module equal to  the cohomology of the 
complex 
\begin{equation}C^*_{\hbar}=\oplus_k 
C^k_{\hbar},\,\,\,C^k_{\hbar}=\oplus_{q-p=k}
\Om^{0,q}(X,\Lambda^pT_X)\widehat\otimes
\C(\hbar),\,\,\db+\hbar\dl:C^k_{\hbar}\to
C^{k+1}_{\hbar} \label{eq:C(h)}
\end{equation}
\kendremark

\subsubsection*{Singularity at $\hbar=0$ of the family of affine structures 
on $\Mh$}\label{s:affstMh}
In this subsection we show that the family of affine 
structures on $\Mh$ have the singularity (\ref{eq:G(h)})
at $\hbar=0$.

 
\begin{prp}\label{pr:pole}
The family of affine structures on 
$\Mh,\hbar\neq 0$ has
pole of order one at $\hbar=0$
\end{prp}
\proof 
A choice of (a class of gauge equivalence of) versal fibered
solution to 
Maurer-Cartan equation $\wgh\in({\V}\,\widehat\otimes 
t_{{\bH}}\C[[t_{\widehat{\bH}}]])^1$  
corresponds 
geometrically to a fiberwise choice of 
coordinates $\rho:\widetilde{{\M}}\to \bH\times \Ac$. 
The affine structure 
on the fiber $\Mh,\hbar\neq 0$ is induced  via period 
mapping to the cohomology
of the complex (\ref{eq:cmplx})
\begin{equation}(t,\hbar)\,\,\,\to \bigl[
\mbox{exp}({i_\hbar(\wgh))}\bigr]\in 
H\otimes(\C(\hbar)\widehat\otimes_{(0)}\Ct)
\label{eq:exp(iwgh)}
\end{equation}
The relative tangent sheaf of the total space of the 
family $p:\MM\to\Ac$
is identified naturally with 
the cohomology of the complex
$(\V\,\widehat\otimes\C[[t_{\widetilde{\bH 
}}]][1],\db+\hbar\dl+[\wgh,\,\cdot\,])$. Denote this 
$\C[[t_{\widetilde{\bH}}]]-$module by ${\mathcal T}_{\MM/\Ac}$.
The Jacobian  of map (\ref{eq:exp(iwgh)})  is equal to
\begin{equation}
\delta\ga \in {\mathcal T}_{\MM/\Ac}\,\,\,\to 
\bigl[i_\hbar(\delta\ga)\wedge
\mbox{exp}(i_\hbar\wgh)\bigr]\in 
H\otimes(\C(\hbar)\widehat\otimes_{(0)}\Ct)
\label{eq:Jac}
\end{equation}
Let us  calculate the affine connection 
\begin{equation}
\Gamma_{ij}^k(\hbar)=
\sum_{\alpha}(\frac{\p a^{-1}}{\p t})^k_\alpha\p_i\p_j a^{\alpha}
\end{equation}
where $\{a^{\alpha}\}$ is a set of linear coordinates on $H$.
 Let  $\p_i\wg(t,\hbar),\,\,
 \p_j\wg(t,\hbar)\in \mathcal 
T_{\MM/\Ac}$  are the cohomology classes 
corresponding to the coordinate 
vector fields 
$\p_i,\,\,\p_j$ tangent to the fibers. 
 Differentiating twice  the expression
(\ref{eq:exp(iwgh)})
 and taking the inverse to the map (\ref{eq:Jac}) we see that
\begin{equation}
{\nabla (\hbar)}_{\p_i}\p_j=\bigl[\frac{1}{\hbar} 
\p_i\wgh\wedge\p_j\wgh +\p_i\p_j\wgh\bigr]\in 
{\mathcal T}_{\MM/\Ac}\otimes\C[\hbar^{-1}]
\label{eq:nbl_ij}
\end{equation}
Family of vector fields 
 $[\p_i\wgh(t,\hbar)\wedge\p_j\wgh(t,\hbar) 
+\hbar\p_i\p_j\wgh(t,\hbar)]$ is regular 
at $\hbar=0$ and does not vanish at $\hbar=0$. Therefore
the order of pole at $\hbar=0$ of the connection 
$\Gamma_{ij}^k(\hbar)$ is 
equal 
to one.
\pendproof 

\begin{prp}\label{pr:res}
  The 
algebra structure defined on $T_{{\M}}$ by the ``residue'' 
of the family of affine connections on $\Mh$ coincides
with the algebra structure defined in \S \ref{ss:product}
\end{prp}
\proof
It follows from the 
equations (\ref{eq:nbl_ij}) and (\ref{eq:product}) 
\pendproof

Conjecturally the family of moduli spaces $\Mh$
should be possible to describe
entirely in terms of the $A_\infty-$ category $D^bCoh(X)$.
The fibers $\Mh$ for $\hbar\neq 0$ equipped with their affine
structure should be identified  with open subsets in periodic
cyclic cohomology of $D^bCoh(X)$. The cyclic cohomology
 corresponds to infinitesimal deformations of a section of $\MM\to\Ac$.

\subsection{Flags and  flat connections over $\AC$}
\label{s:flags}

In this subsection we recall the correspondence due to Rees
(see also \cite{s} )
between  filtrations on a
$\C-$vector space and $\Gm-$equivariant vector bundles on $\AC$.
(Recall that $\Gm$ denotes here the algebraic group whose set
of $\C-$points coincide with $\C\setminus 0$)

\begin{prp}\label{pr:filtr}
The category of finite-dimensional $\C-$vector
spaces equipped with decreasing filtration $(F^{\ge k})_{k\in\Z}$ 
is canonically equivalent to the category of coherent locally free 
$\Gm-$equivariant sheaves over $\AC$.\end{prp}
\proof
Let ${\mathcal V}\to \AC$ be a coherent locally free sheaf equipped with 
an action of $\Gm$ which covers the natural
action of $\Gm$ on the affine line $\AC$.
In particular $\Gm$ acts on the fiber $\CV_0$. This action  defines 
canonical grading 
$\CV_0=\oplus_{i\in \Z} \CV^i_0$ 
where $\Gm$ acts on the graded component $\CV_0^i$ via 
\begin{equation}v\to \lambda^{i}v\,\,\,\,\mbox{for}\,\,\,\lambda\in 
\C^*,\,v\in \CV_0^i\end{equation}
Similarly, any fiber $\CV_{\hbar_0}, \,\hbar_0\neq 0$ receives 
canonical decreasing filtration 
\begin{equation}F^{\geq r}\subset F^{\geq r-1}\dots\subset 
F^{\ge s}=\CV_{\hbar_0},\,r,s\in\Z\end{equation} 
This is the filtration  by  the orders of 
growth at $0\in\AC$ of $\Gm-$invariant sections. Given a choice 
 of isomorphism of ${\mathcal O}_{\AC}-$modules $\CV\simeq\CV_0\otimes
{\mathcal O}_{\AC}$  a $\Gm-$invariant section $v(\hbar)$ can be written 
\begin{equation}v(\hbar)=
\hbar^{-p}v^{(-p)}+\hbar^{-p+1}v^{(-p+1)}+\dots,\,\,v^{(i)}\in 
\CV_0\end{equation} 
where the order of growth $\hbar^{-p}$ and the ``residue'' 
$v^{(-p)}\in \CV_0^{-p}$ does not depend on the choice of the 
isomorphism  $\CV\simeq\CV_0\otimes
{\mathcal O}_{\AC}$.
Put $v(\hbar_0)\in F^{-p}\subset \CV_{\hbar_0}$ if the order at $\hbar=0$ of 
the $\Gm-$invariant section $v(\hbar)$ is less or equal to $\hbar^{-p}$.
The action of $\Gm$ identifies  all the fibers 
$\CV_{\hbar}\,\,\mbox{for}\,\,\hbar\neq0$. This identification respects the 
filtrations. The corresponding associated graded quotent is canonically 
isomorphic to $\CV_0$. 

Conversely, given a vector space 
$\CV_1$ equipped with filtration $(F^{\ge k})_{k\in \Z}$ 
define ${\mathcal 
O}_{\AC}-$module
to be the ${\mathcal O}_{\AC}-$submodule of $\CV_1\otimes 
\C[\hbar^{-1},\hbar]$generated by the elements of the form
$\hbar^{p}v$ where $v\in F^{-p},\,\,p\in\Z$. This is a free finitely 
generated ${\mathcal O}_{\AC}-$module. It has a natural $\Gm-$action 
\begin{equation}\hbar^iv\to\lambda^{i}\hbar^iv,\,\,\mbox{for}\,\lambda\in 
\C^*\end{equation}\pendproof 

More generally, a vector space equipped with filtration having arbitrary
rational
indices $\alpha_i\in \Q$ corresponds  to 
a vector bundle over $\AC$ endowed with a 
  connection over
$\AC\setminus 0$ with regular singularities at zero and at infinity, which
has  monodromy of finite order.

A similar correspondence exists between increasing filtrations and local
systems over $\CP\setminus 0$ having
a pole of the first order at infinity. 
The subspace having an index $\alpha$ 
corresponds to the $\Gm-$invariant sections of order $\leq \hbar^\alpha$
at infinity. 
 
For  projective complex manifold $X^n$ one has the associated Hodge
filtration (\ref{eq:Fgeqr}).
It corresponds to vector bundle
of graded vector spaces ${\mathcal L}^{Hodge}$ over $\AC$
equipped with  connection having regular singularity at $\hbar=0$.
The fiber ${\mathcal L}^{Hodge}_\hbar$ over $\hbar$ is the cohomology of the 
operator $\db+\hbar\p$ acting on
the graded vector space of differential forms. The flat 
multivalued sections
of ${\mathcal L}^{Hodge}$ are 
\begin{equation}
[\phi^{r,q}\hbar^{\frac{r-q+n}{2}}]
\fr \phi^{r,q}\in \Om^{r,q}(X) 
\end{equation}

\subsection{The bundle $\widehat{{\M}^W}\to\CP$ and 
 the map $\Phi^W$.}\label{s:bundle}
In this subsection we give a geometric interpretation
for the results concerning  the generalized
period map  from \S \ref{ch:genper}.

An increasing filtration $W$ 
on $\oplus_i H^i(X,\C)[-i]$ introduced in
\S \ref{s:thegpm} defines 
 an extension of  ${\mathcal L}^{Hodge}$ to a vector bundle
over $\CP$ as we saw it in \S \ref{s:flags}. 
\begin{prp}\label{pr:compl}
A filtration $W$ complementary to the Hodge filtration  in the sense of
eq.~(\ref{eq:compl})
 defines a trivial sheaf over $\CP$.
\end{prp}
\proof  It follows from (\ref{eq:compl})
 that $\oplus_i H^i(X,\C)[-i]=\oplus_s(F^{\ge s}\cap W_{\le s})$
The trivialization of the extension of the sheaf ${\mathcal L}^{Hodge}$ is 
given by the  sections $\hbar^{-s} v(\hbar)$ where $v(\hbar)\in F^{\ge s}\cap 
W_{\le s}$
 is a flat (multivalued) section.
\pendproof 

\subsubsection*{The extension of the bundle ${\MM}\to\AC$.}\label{ss:exten}
We explain here that the filtration $W$ defines canonical extension 
$\widehat {{\M}^W}\to\CP$ of the bundle $\MM\to\AC$.

Let $\pi_\hbar$ denotes the base point of the formal moduli space $\Mh$.
The $\Z-$graded tangent space $T_{\pi_\hbar}{\Mh}$ at the base point is 
identified naturally with the cohomology of the complex 
$(\V,\db+\hbar\dl)[1]$.
A choice of nowhere vanishing holomorphic $n-$form 
 induces  an isomorphism 
of $\Z-$graded vector spaces  $T_{\pi_\hbar}{\Mh}\simeq
{\mathcal L}^{Hodge}_\hbar[n]$. The isomorphism is defined canonically
up to a multiplication by nonzero constant. 
Therefore we get the structure of a local system with regular
singularity at $\hbar=0$ on the bundle 
$T_{\pi_\hbar}\to\AC$.
Let us denote by ${\mathcal D}_{\p/\p\hbar}^{0}$ the 
corresponding
 connection.
The flat multivalued sections of the local system can be written as
\begin{equation}[\sum_{p,q}\hbar^{-\frac{p+q}{2}+n}\ga^{p,q}]
\fr \,\sum_{p,q}\ga^{p,q}\in\Ker\,\db+\dl,\,\,\ga^{p,q}\in 
\Om^{0,q}(X,\Lambda^pT_X)
\end{equation}
The filtration $F$ maps to the filtration 
$$(F^H)^{\ge r}=\{[\ga]|\ga\in\oplus_{p+q\leq 
2(n-r)}\Om^{0,q}(X,\Lambda^pT_X)\}$$ on the graded vector space 
$H=T_{\pi_1}{\M}_1$. The filtration $W$ maps similarly to 
a filtration 
$W^H$  and
defines an 
extension $\widehat{T^W_{\pi_\hbar}}\to\CP$ of the bundle 
$T_{\pi_\hbar}\to\AC$.
The affine structure on the fibers $\Mh,\,\hbar\neq 0$ gives an 
identification of $\Mh\fr \hbar\neq 0$ with formal neighborhood of zero in
$T_{\pi_\hbar}\Mh$. Therefore the filtration $W$ induces  extension
$\widehat{{\M}^W}\to\CP$ of the bundle of formal manifolds $\MM\to\AC$.
Notice that the family of affine structures on $\Mh$ continues analytically
at $\hbar=\infty\in \CP$.

Suppose that we are given a  complex analytic bundle over 
$\CP$ and a  section with analytically trivial normal bundle.
Then it is easy to see that the formal neighborhood
of this section is trivial as a nonlinear bundle and the 
trivialization is unique. It follows from the fact that
$H^1(\CP,G)$=0 where $G$ is the projective limit
of  nilpotent algebraic groups of jets of formal diffeomorphisms
of a superspace with the first jet the same as of the identity 
diffeomorphism.
The fact about  $H^1$ follows from 
$H^1(\CP,\mathcal{O})=H^2(\CP,\mathcal{O})=0$. 

The theorem \ref{th:mainth} and the proposition \ref{pr:well-def} 
describe the analog of this phenomena in the
category of formal bundles over $\CP$.
We saw above that a choice of
fibered versal solution $\wgh$ corresponds
to a fiberwise choice of coordinates on $\MM\to\AC$ and
the formal power series $\exp(l_\hbar\wgh)$ represents 
the map to the affine space associated with the
affine structure on the fibers $\Mh$. The theorem
\ref{th:mainth}
is interpreted then as the 
fact 
 that the nonlinear formal bundle 
$\widehat{\M^W}$ is trivial.

\subsection{Canonical connection on $\MM\to\AC$ and symmetry vector field.}

In this subsection we will see that there exists a natural 
symmetry vector field acting on all the structures on ${\M}$.
It is given by the ``residue'' of the canonical connection on the
family $\Mh$.

Let us consider the 
following connection on $T_{\pi_\hbar}$:
\begin{equation}
\mathcal{D}^T_{\p/\p\hbar}=\mathcal{D}^{0}_{\p/\p\hbar}+
\frac{(n-1)}{\hbar}
\end{equation}
The flat multivalued sections of this local system
are
\begin{equation}[\sum_{p,q}\hbar^{-\frac{p+q}{2}+1}\ga^{p,q}]
\fr \,\sum_{p,q}\ga^{p,q}\in\Ker\,\db+\dl,\,\,\ga^{p,q}\in 
\Om^{0,q}(X,\Lambda^pT_X)\label{eq:locst} 
\end{equation}
This local system corresponds to the
same filtration $F^H$ having all the indexes shifted by $(1-n)$:
\begin{equation}
(F^H)^{\geq 1}\subset (F^H)^{\geq \frac{1}{2}}\subset \dots 
\subset (F^H)^{\geq 1-n}
\end{equation}

Let ${\mathcal D}_{\p/\p\hbar}$ denotes the connection
on the bundle $p:\MM\to\AC$ 
induced via the canonical map $\Mh\to T_{\pi_\hbar}$
from the connection ${\mathcal D}^T_{\p/\p\hbar}$.
In other words by definition the covariant derivative
along $\p/\p\hbar$ is the following 
formal series  
\begin{equation} {\mathcal D}_{\p/\p\hbar}:=\frac{\p}{\p\hbar}-
\sum_{a,b}\frac{\p (\Phi^H)^{a}}{\p \hbar}
\bigr(\frac{\p\Phi^{-1}}{\p\tau}\bigl)^{b}_a
\frac{\p}{\tau^b} \end{equation}
where 
$\Phi^H(\tau,\hbar)=[l_\hbar(\hbar\exp
(\frac{1}{\hbar}\wgh)-\hbar)]=[\mbox{exp}(l_\hbar\wgh)-1]$ 
is the map $\Mh\to T_{\pi_\hbar}$ written in the coordinates
corresponding to a fibered
versal solution $\wgh$ and the frame on $T_{\pi_\hbar}$
which is covariantly constant with respect to
$\mathcal{D}_{\p/\p\hbar}^T$.

\begin{prp}\label{pr:E(h)}
The connection $\mathcal{D}_{\p/\p\hbar}$
has regular singularity at $\hbar=0$, in other words
\begin{equation}\Dh=
\frac{\p}{\p\hbar}+\hbar^{-1}E_{-1}(\tau)^\beta\frac{\p}{\tau^\beta}
+E_{0}^\beta(\tau)\frac{\p}{\tau^\beta}+\hbar E_{1}+\ldots 
\label{eq:E(h)}\end{equation}
\end{prp}
\proof
The proof is parallel to the proof of the prop.~\ref{pr:pole}.
The tangent space map induced by $\Phi$ is
\begin{equation}
\frac{\p\Phi}{\p\tau}:\delta\wgh\to 
l_\hbar\delta\wgh\,e^{l_\hbar\wgh}
\end{equation}
The image of $\p\Phi/\p\hbar$ under the inverse map is
\begin{equation}
{\frac{\p\Phi}{\p\tau}}^{-1}\frac{\p\Phi}{\p\hbar}=l_{\hbar}^{-1}
\p_\hbar (l_\hbar) \wgh +\p_\hbar\wgh
\in  
{\mathcal T}_{\MM/\Ac}\otimes\C[\hbar^{-1}] 
\end{equation}
Notice that
\begin{equation}l_{\hbar}^{-1}
\p_\hbar 
(l_\hbar):\,\,\,[\sum_{p,q}\ga^{p,q}]\to\sum_{p,q}{\frac{p+q-2}{2\hbar}}
[\ga^{p,q}] \end{equation}
Therefore the connection $\Dh$ takes the form (\ref{eq:E(h)}) with
\begin{equation}
E_{-1}(\tau)=[\sum_{p,q}\frac{2-p-q}{2}
\wgh^{p,q}|_{\hbar=0}(\tau)]\in 
{\mathcal T}_{{\M}}
\end{equation}
where ${\mathcal T}_{{\M}}$ is identified
with the cohomology of the complex
$(\V\wCtW[1],\db+[\wg|_{\hbar=0},\cdot])$.
\pendproof 

It follows easily from the
transformation law of the connection $\Dh$ under the fiberwise coordinate
changes that the ``residue'' part $E_{-1}(\tau)$ is a well defined
vector field on ${\M}$. In the sequel we will often denote
this vector field
simply by $E(\tau)$.

\begin{prp}
The vector field $E(\tau)$ acts
as a conformal symmetry
on the algebra structure on ${\mathcal T}_{{\M}}$:
 ${Lie}_E(\circ)=\circ$, that is, for any vector fields 
$u,v\in{\mathcal T}_{{\M}}$ 
\begin{equation}
[E,u\circ v]-[E,u]\circ v-u\circ[E,v]=u\circ v\label{eq:Ecirc} 
\end{equation}
\end{prp}
\proof 
The connection $\Dh$ respects the affine structure on the fibers 
$\Mh$.
 The compatibility of ${\mathcal D}_{\p/\p\hbar}$  with the 
 affine structure 
can be written as 
\begin{equation}{\mathcal D}_{\p/\p\hbar}(\nabla_{u_\hbar}(v_\hbar))=
\nabla_{{\mathcal D}_{\p/\p\hbar}(u_\hbar)}(v_\hbar)+
\nabla_{u_\hbar}({\mathcal D}_{\p/\p\hbar}(v_\hbar))\label{eq:Dnbl}
\end{equation}
where  $u_\hbar=u+u_{(1)}\hbar+\ldots,\,\,\,\,v_\hbar=v+v_{(1)}\hbar 
+\ldots\,\,\,\,\in
\G(\AC,p_*{\mathcal T}_{\MM/\AC})$.
 Taking the terms of order $\hbar^{-2}$ at both sides of 
 the formula (\ref{eq:Dnbl})
one obtains the equation (\ref{eq:Ecirc}).
\pendproof

\begin{prp}\label{pr:Ei=0}
Assuming that $\wgh^W$ is the solution which satisfies the conditions 
of the theorem \ref{th:mainth}
one has 
\begin{equation}\Dh=
\frac{\p}{\p\hbar}-\hbar^{-1}s(r_a)
\tau^a\frac{\p}{\p\tau^a}\label{eq:Ei=0}
\end{equation}
where $s(r_a)=-r_a+n-1$ for $\tau^a\in
(W_{\le r_a}/W_{\le r_a-\frac{1}{2}})^{dual}$.
In particular, one has
 $E_0=E_1=\dots=0$ in the coordinates corresponding to 
$\wgh^W$.\end{prp}
\proof
Let us  notice that
\begin{gather}
\bigl({\frac{\p\Phi}{\p\tau}}^{-1}\bigr)^a_{b}
\in\hbar^{-s(r_b)}\C[\hbar^{-1}]
\wCtW\notag\\
\frac{\p\Phi^{b}}{\p\hbar}\in\hbar^{s(r_b)-1}\C[\hbar^{-1}]
\wCtW\notag
\end{gather}
Then the  proposition \ref{pr:E(h)}
implies that
\begin{multline}
E(\hbar)=-\sum_{a,b}
\bigl(\frac{\p(\Phi^{b}_{(0)}\hbar^{s(r_b)})}
{\p\tau^a}\bigr)^{-1}
\frac{\p(\Phi^{b}_{(0)}\hbar^{s(r_b)})}{\p\hbar}\frac{\p}
{\p\tau^a}=\\
=-\hbar^{-1}\sum_{a} s(r_a)
\tau^a\frac{\p}{\p\tau^a}
\end{multline}
\pendproof 

\subsection{Invariants of the 
generalized VHS.}\label{s:invVHS}

Fix a point
$[X_z]$ in the classical moduli space of complex structures
on $X_{C^\infty}$ and
an increasing  filtration $W$ on $\oplus_k H^k(X_z,\C)[-k]$
which is complementary to the Hodge filtration. 
Then the coefficients in the Taylor expansion
of the power series 
representing the $3-$tensor  $A^{a}_{bc}(\tau_W)$ written in coordinates
$\{\tau_W\}$
form collection of polylinear $S_k-$symmetric maps
\begin{equation}
A^{(k)}_z:\,\,\,{\mathcal{H}}^{\otimes k}\,\,\to\,\,
\mathcal{H},\,\,\,
\,\,\, 
\mathcal{H}=\oplus_{p,r}H^r(X,\Lambda^{p}T_X)[p-r]\label{eq:invVHS}
\end{equation}
where we have used the natural identification of graded vector spaces
$\mathcal{H}
\simeq\mbox{Gr}\,W^H\stackrel{\vdash\Om_{X}}{\simeq}\mbox{Gr}\, W$

In the view of the proposition \ref{pr:filtr} it is
natural to consider the family of moduli spaces 
$\Mh$ equipped with the canonical connection
$\mathcal{D}_{\p/\p\hbar}$ as a nonlinear analog of 
the Hodge filtration. On the other hand as we will explain
below such families of the moduli spaces can be canonically
identified when the base point $\pi\in\M^{{classical}}$ changes.
It is natural to consider the family of maps (\ref{eq:invVHS})
as invariants of the generalized variations of Hodge structure.
We will see that these invariants form  exactly the same structure as
the Gromov-Witten invariants.

\subsubsection*{Potentiality of $A_{ab}^{c}(\tau_{W})$}
Here we start to study the properties
of the invariants of generalized VHS.   
As one of the consequences we show that
the 
tensor of algebra structure constants $A^c_{ab}$ written using
the  coordinates $\{\tau_W\}$ on ${\M}$  satisfies
\begin{equation} \forall\,\, 
\,\,a,b,c,d\,\,\,\,\,\,\,\,\p_aA^b_{cd}=(-1)^{\bar a\bar c}\p_cA^b_{ad}
\label{eq:dA=0} 
\end{equation}

\begin{thm}\label{thm:dA=0}
The power series  representing the family of affine connections on 
the moduli spaces $\Mh$ is written using the 
fiberwise coordinates corresponding to $\wg^W(\tau,\hbar)$ as
\begin{equation}\nabla(\hbar)=d+\hbar^{-1}A_{ab}^c(\tau) 
\end{equation}\end{thm}
\proof
Let 
\begin{gather}
\begin{split}
&(\Phi^W)^{a} 
(\tau_{W})=\bigr[e^{
l_\hbar\wg^W(\tau,\hbar)}-1\bigl]^{a}=\\ 
&\Phi^{a}_{(0)}(\tau)\hbar^{s(r_a)}+\Phi_{(-1)}^{a}
(\tau)\hbar^{s(r_a)-1}+\dots \label{eq:phiWa}
\end{split}\\
\mbox{where}\,\,\,\,(\Phi^W)^{a}_{(-i)}\in  
\goth M^{i+1}_{\C[[{\tau}_{W}]]}\,\,\,\,\mbox{and}\,\,\, 
(\Phi^W)_{(0)}^{a}(\tau)=\tau^{a}\notag
\end{gather}
 are the components of the map (\ref{eq:phiW}). 
 In particular 
\begin{equation}
\frac{\p}{\tau^b}\frac{\p}{\tau^c}\Phi_{(0)}^{a}(\tau)=0
\end{equation} 
for any $a,b,c$. It follows from the equation (\ref{eq:phiWa}) that
\begin{equation}\bigr({\frac{\p\Phi} 
{\p\tau}}^{-1}\bigl)^{b}_{a}\in\hbar^{-s(r_a)}\C[[\hbar^{-1}]] 
\end{equation}
Therefore 
\begin{equation}
\G_{bc}^e(\hbar)=\sum_{a}\p_b
\p_c\Phi^{a}\bigr({\frac{\p\Phi} 
{\p\tau}}^{-1}\bigl)^{e}_{a}
\in\hbar^{-1}\C[[\hbar^{-1}]]\wCtW 
\end{equation}
On the other hand according to the proposition \ref{pr:pole}
\begin{equation}\G_{bc}^e(\hbar)
\in \hbar^{-1}\C[[\hbar]]\widehat\otimes\C[[\tau_{W}]]\end{equation}
\pendproof 

\begin{cor}
The tensor of algebra structure constants $A_{ab}^c(\tau_W)$ written
in the coordinates on ${\M}$ corresponding to the solution 
$\wg^W|_{\hbar=0}$
has the property (\ref{eq:dA=0}).
\end{cor}
\proof
The tensor $A_{ab}^c(\tau)$ can be interpreted as the ``residue'' of 
the affine connection on $\Mh$
(see \ref{pr:res}).
Consider the equation (\ref{eq:dG=0})
written in the fiberwise coordinates corresponding to the 
solution $\wg^W(\tau,\hbar)$.
This equation expresses the fact that $\nabla(\hbar)$ is flat. 
The term of order $\hbar^{-1}$ gives 
\begin{equation}d(\sum_a A_{ab}^c d\tau^a)=0
\end{equation}
which is exactly the required property (\ref{eq:dA=0}).
\pendproof

\subsection{Flat metric on ${\M}$}

Here we explain that if the filtration $W$ on $H^*(X,\C)$ is isotropic with 
respect to the Poincare pairing  then a (holomorphic)
flat metric which  is  compatible with the
algebra structure (equation (\ref{eq:g(uvw)})) and the affine structure 
$\{\tau_W\}$
(equation (\ref{eq:dg=0})) is induced on ${\M}$.

Assume that the filtration $W_{\le 0}\subset 
W_{\le \frac{1}{2}}\subset\dots\subset W_{\le n}=\oplus_k H^k(X,\C)[-k]$ 
introduced in \S \ref{ss:normsol} is  isotropic
with respect to the Poincare pairing $<\,,\,>$ in the following sense:
\begin{equation}
<\phi_1,\phi_2>=0\fr\,\phi_1\in W_{\le r},\,\phi_2\in 
W_{\le n-r-\frac{1}{2}} \label{eq:isotr}
\end{equation}
This property together with the corresponding property 
\begin{equation}
<\phi_1,\phi_2>=0\fr\,\phi_1\in F^{\ge r},\,\phi_2\in 
F^{\ge n-r+\frac{1}{2}} \end{equation}
of the
Hodge filtration implies
\begin{equation}
<\phi_s,\phi_{s'}>\neq 0\fr \, \phi_r\in F^{\ge r}\cap 
W_{\le r},\,\phi_{r'}\in F^{\ge r'}\cap W_{\le r'}\Longrightarrow 
r+r'=n \label{eq:orth}\end{equation}


Assume that a choice of nowhere vanishing holomorphic $n-$form $\Om_{X}$
is fixed.
Define the pairing
\begin{equation}
(\ga_1,\ga_2):=\int_X(\ga_1\wedge\ga_2)\vdash\Om_{X}\wedge\Om_{X}\fr\,\,
\ga_1,\,\ga_2\in \oplus_{p,q}\Om^{0,q}(X,\Lambda^pT_X)\label{eq:pair} 
\end{equation}

\begin{prp}
\begin{gather}
(\db\ga_1,\ga_2)=(-1)^{p_1+q_1+1}(\ga_1,\db\ga_2)\\
(\dl\ga_1,\ga_2)=(-1)^{p_1+q_1}(\ga_1,\dl\ga_2)
\end{gather}
\end{prp}
\proof
The proof follows from the integration by parts formulas.
\pendproof 

\begin{prp}\label{pr:ThT-h}
The pairing (\ref{eq:pair})
 induces  a natural pairing on the opposite fibers of the 
local system $T_{\pi_\hbar}\otimes T_{\pi_{-\hbar}}\to\C$
\end{prp}
\proof
Recall that
 the fibers of the local system $T_{\pi_\hbar}\to\AC$ are naturally 
identified with the cohomology of the complex $(\V,\db+\hbar\dl)[1]$.
\pendproof 

\begin{prp}
The pairing $\hbar^{n-2}(\,,\,):
T_{\pi_\hbar}\otimes T_{\pi_{-\hbar}}\to\C$
 is locally constant:
\begin{equation}\p_\hbar\hbar^{n-2}(u_\hbar,v_{-\hbar})=
\hbar^{n-2}({\mathcal D}_{\p/\p\hbar}^T u_\hbar,v_{-\hbar})+
\hbar^{n-2}(u_\hbar,{\mathcal D}^T_{-\p/\p\hbar}v_{-\hbar}) \end{equation}
where  $u_\hbar=u+u_{(1)}\hbar+\ldots,\,\,\,\,v_\hbar=v+v_{(1)}\hbar 
+\ldots\,\,\,\,\,\,\in \G(\AC,\pi_\hbar^*T_{\MM/\AC})$. 
\end{prp}
\proof
The locally constant elements of 
$\G(\AC\setminus 0,\pi_\hbar^*T_{\MM/\AC})$ are of the form 
(\ref{eq:locst}).
\pendproof 
Consider the same pairing written in the locally constant frame
on $T_{\pi_\hbar}\to\AC$.
\begin{prp} \label{pr:(u,v)infty}
Assuming that the filtration $W$ is isotropic one has
\begin{equation}\forall \,\,\,u,v\,\,\,\in \oplus_r 
W_{\le 
r}^H\hbar^{s(r)}\C[[\hbar^{-1}]]\,\,\,\,\,\,\,\,\,(u(\hbar),v({-\hbar}))\in 
\C[[\hbar^{-1}]] \end{equation}In other words the pairing 
$(\,,\,):T_{\pi_\hbar}\otimes 
T_{\pi_{-\hbar}}\to\C$ extends over 
$\hbar=\infty$. 
\end{prp}
\proof
For  $u\in H$ the corresponding element of $(\db+\hbar\dl)-$cohomology
class can be written as $l_{\hbar}^{-1}u$, where 
$l_\hbar^{-1}:[\sum\ga^{p,q}]\to[\sum\hbar^{-\frac{p+q-2}{2}}\ga^{p,q}]$.
Notice that for any $u,v\in H$
\begin{multline}
\int_X(l_\hbar^{-1}u\wedge 
l_{-\hbar}^{-1}v)\vdash\Om\wedge\Om=\\
=(-1)^{1+\overline{v}(n-\frac{1}{2})}\int_X (l_\hbar^{-1}u\vdash\Om)\wedge
 (l_\hbar^{-1}v\vdash\Om)\\
=(-1)^{1+\overline{v}(n-\frac{1}{2})}
\hbar^{2-n}\int_X (u\vdash\Om)\wedge(v\vdash\Om)\label{eq:uomvom}
\end{multline}
where $\mbox{deg}\, v^{p,q}=q-p$ for $v^{p,q}\in \Om^{0,q}(X,\Lambda^p T_X)$.
To complete the proof it is enough to notice that 
by the property (\ref{eq:orth}) for $u\in W_{\le r},v\in W_{\le r'}$
\begin{equation}\int_X (u\vdash\Om)\wedge(v\vdash\Om)\neq 0
\Longrightarrow s(r)+s(r')
\leq n-2\end{equation}
\pendproof

Let us denote by 
\begin{equation}\Phi^{W,H}(\tau_{W},\hbar)
=\bigr[\hbar\mbox{exp}\,
{\frac{1}{\hbar}\wg^W(\tau,\hbar)}-\hbar\bigl]\in H_{(\hbar)}
\widehat{\otimes}_{(0)}\C[[\tau_W]] 
\end{equation}
 the formal power series representing the 
canonical map $\Mh\to T_{\pi_\hbar}$ written 
using the fiberwise coordinates corresponding to $\wgh^W$. 

Consider the pairing induced on the tangent sheaf of ${\M}$ via  this
 map 
from the pairing $T_{\pi_\hbar}\otimes T_{\pi_{-\hbar}}\to\C$
\begin{equation}(\Phi^W)^*(u,v):=
(\p_u\Phi^{W,H}(\tau_{W},\hbar),
\p_v\Phi^{W,H}(\tau_{W},-\hbar))
\label{eq:def}
\end{equation}
\begin{thm}\label{th:(,)inh0}
Assuming that the filtration $W$ is isotropic
the induced pairing does not depend on $\hbar$:
\begin{equation}\forall \,\,\,u,v\in 
{\mathcal T}_{{\M}_0}
\,\,\,\,\,\,(\Phi^W)^*(u,v)\in \hbar^0\C[[\tau_{W}]] \end{equation}
\end{thm}

\proof
Notice that 
\begin{gather}
(\Phi^W)^*(u,v)=
\int_X 
\bigr(\p_u\wg^W(\tau,\hbar) e^{\frac{1}{\hbar}\wg^W(\tau,\hbar)}\wedge
\p_v\wg^W(\tau,-\hbar) e^{-\frac{1}{\hbar}
\wg^W(\tau,-\hbar)}\bigl)\vdash\Om\wedge\Om=\\
\int_X
\bigr(\p_u\wg^W(\tau,\hbar) \p_v\wg^W(\tau,-\hbar)
e^{\frac{1}{\hbar}((\widetilde{\ga}^W_{(0)}(\tau)+\hbar
\widetilde{\ga}^W_{(1)}(\tau)+\ldots)-
(\widetilde{\ga}^W_{(0)}(\tau)-\hbar\widetilde{\ga}^W_{(1)}(\tau)+\ldots))
}\bigl)\vdash\Om\wedge\Om\label{eq:(,)inh0}
\end{gather}
Hence $(\Phi^W)^*(u,v)\in \C[[\hbar]]\wCtW$

On the other hand 
for any $u\in 
{\mathcal T}_{{\M}_0}$
its  image under the map $(\Phi^{W,H})_*$ written 
in the locally constant frame on 
$T_{\pi_\hbar}\to\AC$ lies in
\begin{equation}\Phi^W_*u\in \oplus_r 
W^H_r\hbar^{s(r)}\C[[\hbar^{-1}]]\wCtW\end{equation}
Therefore 
by the proposition \ref{pr:(u,v)infty}
\begin{equation}\forall u,v\in {
\mathcal T}_{{\M}_0}\,\,\,\,\,\,
(\Phi^W)^*(u,v)\in 
\C[[\hbar^{-1}]]\widehat\otimes\C[[\tau_{W}]] 
\end{equation}
\pendproof

\begin{prp}\label{pr:g(uvw)}
The pairing $(\Phi^W)^*:T_{{\M}_0}^{\otimes 2}\to\C$ is 
compatible with the algebra structure on $T_{{\M}_0}$ 
 in the following sense:
\begin{equation}\forall \,\,\,u,v,w
\in {\mathcal T}_{{\M}_0}\,\,\,\,\,\,\,\,\,(\Phi^W)^*(u\circ 
w,v)=(\Phi^W)^*(u,w\circ v) \label{eq:g(uvw)}
\end{equation}
\end{prp}
\proof 
The pairing $T_{\pi_\hbar}\otimes T_{\pi_{-\hbar}}\to\C$ induces
via the canonical map $\Mh\to T_{\pi_\hbar}$
the pairing \begin{equation}\Phi_\hbar^*:{\mathcal T}_{\Mh}\otimes 
{\mathcal T}_{{\M}_{-\hbar}}\to{\mathcal O}_{\Mh\times{\M}_{-\hbar}} 
\label{eq:phih}
\end{equation} 
One can regard the pairing 
$(\Phi^W)^*:{\mathcal T}_{{\M}_0}^{\otimes 2}\to{\mathcal O}_{{\M}_0}$
 as being induced
from the pairing $\Phi_\hbar^*$ via the identification of the
fibers $\Mh\sim{\M}_0$ associated with the fiberwise coordinate
choice corresponding to $\wgh^W$.
The pairing $\Phi_\hbar^* (\,,\,)$
is compatible with the affine connections on the fibers 
$\Mh,{\M}_{-\hbar}$:
\begin{multline}
\forall \,\,\,\,u,v,w\in \,\,\,\G(\MM,{\mathcal T}_{\MM/\AC})
\label{eq:nbl(,)}\\
Lie_{ 
 (w(\hbar),w(-\hbar))} \Phi_\hbar^*(u,v)
 =\Phi_\hbar^*(\nabla(\hbar)_w u(\hbar), v(-\hbar))+ 
 \Phi_\hbar^*(u,\nabla(-\hbar)_w 
v(-\hbar))
\end{multline}
To prove this equality it is enough to notice that
it becomes trivial if one uses affine coordinates on the
fibers $\Mh,{\M}_{-\hbar}$.
For $u,v,w\in  \G({\M}_0,{\mathcal T}_{{\M}_0})$ 
denote by $u(\hbar), v(\hbar), w(\hbar)\in 
\G(\MM,{\mathcal T}_{\MM/\AC})$ the relative vector fields   such 
that $u(0)=u,v(0)=v,w(0)=w$ and $u(\hbar), v(\hbar), w(\hbar)$ are constant 
in the trivialization of $p$ defined by $\wgh^W$.
 The equation (\ref{eq:g(uvw)})
is obtained from (\ref{eq:nbl(,)})
applied to $u(\hbar),v(\hbar),w(\hbar)$ by looking at 
the terms of order 
$\hbar^{-1}$. 
\pendproof

\noindent{\it Second proof.}
Taking the terms of zero order in $\hbar$ in (\ref{eq:(,)inh0}) we see that
\begin{equation}
(\Phi^W)^*(u,v)=\int_X(\p_u\wg^W_{(0)}\wedge\p_v\wg^W_{(0)}\wedge
e^{2\wg^W_{(1)}(\tau)})
\vdash\Om\wedge\Om\end{equation}
Therefore 
\begin{multline}
\forall \,\,\,u,v,w\,\,\in{\mathcal T}_{{\M}_0}\,\,\,\,
(\Phi^W)^*(u\circ w,v)=\\
\int_X(\p_u\wg^W_{(0)}\wedge\p_v\wg^W_{(0)}\wedge\p_w\wg^W_{(0)}\wedge
e^{2\wg^W_{(1)}(\tau)})\vdash\Om\wedge\Om=
(\Phi^W)^*(u,w\circ v)\label{eq:Aabc}
\end{multline}
\pendproof

Recall that $\{\Delta_a\}$ denotes a basis
in the graded vector space $H$ compatible
with the direct sum decomposition 
$H^*(X,\C)=\oplus_r F^{\ge r}\cap W_{\le r}$. We 
have also denoted via $\{\tau^a_{W}\}$  the 
natural coordinates on ${\M}$ assocated with the formal power series 
$\wg^W|_{\hbar=0}(\tau)$. 
Since the filtration $W$ is isotropic 
the Poincare pairing induces the pairing 
on $\mbox{Gr}\,W_*$ which we will denote by the
same symbol $<\,,\,>$.

\begin{prp}
\begin{equation}(\Phi^W)^*
(\p_a,\p_b)=
(-1)^{\bar{b}(n-\frac{1}{2})-r_b+n}<\Delta_a\vdash\Om,\Delta_b\vdash\Om> 
\label{eq:dg=0}\end{equation}
where  $\Delta_b\in W^H_{\le r_b}\cap (F^H)^{\geq 
r_b}$ and $\bar{b}=k$ for $\Delta_b\in H^k$.
In particular, the pairing 
$(\Phi^W)^*:({\mathcal T}_{{\M}_0})^{\otimes 2}\to{\mathcal O}_{{\M}_0}$ is 
constant in the coordinates $\{\tau_{W}\}$. \end{prp}
\proof
It follows from the theorem \ref{th:(,)inh0} that one should
take into account only the following components of $\Phi^W$ in 
(\ref{eq:def}):
\begin{multline}(\Phi^W)^*
(\p_a,\p_b)=\sum_{c,d}\hbar^{s(r_c)}(-\hbar)^{s(r_d)}
([\p_a\Phi^{c}_{(0)}\Delta_c],[\p_b\Phi^{d}_{0}\Delta_d])\\
=(-1)^{\bar{b}(n-\frac{1}{2})-r_b+n}<\Delta_a\vdash\Om,\Delta_b\vdash\Om>
\end{multline}
where in the last step we have used  the property (\ref{eq:PhiW0}) and the 
calculation (\ref{eq:uomvom}).
\pendproof

\begin{cor} 
The pairing $(\Phi^W)^*(\,,\,)$ is symmetric and 
nondegenerate.
\end{cor}
\proof 
These  properties follow from the equation (\ref{eq:g(uvw)})
and (\ref{eq:dg=0}) correspondingly.

\begin{prp}
The vector field $E\in{\mathcal T}_{{\M}_0}$
is conformal with respect to the metric induced by  $\Phi^W$ :
$Lie_ E(\Phi^W)^*(\,,\,)=(2-n)(\Phi^W)^*(\,,\,)$, that is
\begin{multline}
\forall \,\,u,v\in \,\,{\mathcal T}_{{\M}_0}\\
E(\Phi^W)^*(u,v)-(\Phi^W)^*([E,u],v)-(\Phi^W)^*
(u,[E,v])=(2-n)(\Phi^W)^*(u,v)\label{eq:EPhi}
\end{multline}
\end{prp}
\proof
As in the proof of proposition \ref{pr:g(uvw)} let us regard the 
pairing $(\Phi^W)^*:{\mathcal T}_{\M}^{\otimes 2}\to{\mathcal O}_{\M}$ 
as being induced from the pairing $\Phi_\hbar^*$ defined in (\ref{eq:phih}).
The connection $\Dh$
respects the pairing $\hbar^{n-2}\Phi_\hbar^*$:
\begin{multline}\Dh 
\hbar^{n-2}\Phi_{\hbar}^*(u(x,\hbar),v(y,-\hbar))=\\
\hbar^{n-2}\Phi_{\hbar}^*(\Dh 
u(x,\hbar),v(y,-\hbar))+\hbar^{n-2}\Phi_{\hbar}^*(u(x,\hbar),{\mathcal 
D}_{-\p/\p\hbar}v(y,-\hbar))\label{eq:DPhi}\end{multline}
Let us apply this equality to $u(\hbar),v(\hbar)\in \G(\AC,p_*
{\mathcal T}_{\MM/\AC})$
which are constant in the trivialization  corresponding to 
$\wgh^W$. Notice that for  such vector fields 
\begin{equation}(\Phi^W)^*(u(\tau,0),v(\tau,0))=\Phi^*_\hbar
(u(\tau,\hbar),v(\tau,-\hbar)) \end{equation}
Let us restrict both sides of 
the equation (\ref{eq:DPhi}) to $x=y$ and pick up the terms
of order $\hbar^{n-3}$. A short calculation gives 
the equation (\ref{eq:EPhi}).
\pendproof

\subsection{Flat identity}

Consider the coordinates on ${\M}$ associated with the  
power series $\wg^W_{(0)}=\wg^W|_{\hbar=0}(\tau_{ W})$.
\begin{prp}\label{pr:1}
The coordinate vector field $\p_0$ corresponding to the element
$\Delta_0=[1]\in  W^H_{\le n}/W_{\le n-\frac{1}{2}}^H$ 
is the identity with respect to
 the algebra structure on ${\mathcal T}_{{\M}}$.
\end{prp}
\proof
We need to prove that
\begin{equation}
\p_0\wg^W_{(0)}=[1] \,\,\,\in \,\,\,\Ker\, \db+[\wg^W_{(0)},\,]/ 
\mbox{Im}\, \db+[\wg^W_{(0)},\,]  \end{equation}
The
relative
vector field corresponding 
to 
\begin{equation}
\delta\wg=\hbar[1]\in \Ker \,\db+\hbar\dl+[\wg^W,\,]/ \mbox{Im}\, 
\db+\hbar\dl+[\wg^W,\,]
\end{equation}
is sent
by the map 
\begin{equation}
[\exp(l_{\hbar}\wg^W)]_*:{\mathcal T}\Mh\to {\mathcal 
T}T_{\pi_\hbar}\label{eq:exp*}
\end{equation} to the  
family of eulerian vector fields
\begin{equation}
\bigl[\frac{\delta\exp(l_\hbar\wg^W)}{\delta \wg}\bigr]
=[\exp(l_\hbar\wg^W)]=x_\hbar\frac{\p}{\p x_\hbar}
\end{equation}
 in the fibers of the local system
$T_{\pi_\hbar}$.
On the other hand, since the series $\wg^W$ satisfies the
conditions of the theorem \ref{th:mainth}
the image of the family of eulerian vector fields
under the inverse to the map (\ref{eq:exp*})
is the family written in the coordinates $\tau_W$ as
\begin{equation}
\sum_{a,b}(1+\Phi^{W,a})
\bigr(\frac{\p\Phi^{W}}{\p\tau}^{-1}\bigl)^b_a\frac{\p}{\p\tau^b}=
\hbar\frac{\p}{\p\tau^0}+\hbar^{0}v_0(\tau)+
\hbar^{-1}v_{-1}(\tau)+\dots
\end{equation}
Therefore 
\begin{equation}
\frac{\p}{\p\tau^0}=[1]\in 
\Ker \,\db+\hbar\dl+[\wg^W,\,]/ \mbox{Im}\, 
\db+\hbar\dl+[\wg^W,\,]
\end{equation}
\pendproof

\subsection{Dependence on the base point $\pi\in{\mathcal M}^{classical}$}

In this subsection we study the variation of the
structure introduced above with respect to the deformations of the
base point $\pi$.

Let us denote by 
$\phi$ the coordinates defined  by the map (\ref{eq:coor})
on a neighborhood of the point $[X_0]$
of the moduli space $\M^{classical}$
of complex structures on $X_{C^\infty}$.

Let $W(\phi),W(0)=W$ be the family of filtrations
on $\oplus_i 
H^i(X_\phi,\C)$ locally constant  with respect 
to the Gauss-Manin 
connection. If $W(0)$ is complementary
to the Hodge filtration at $\phi=0$ the filtration
$W(\phi)$ remains complementary to the Hodge filtration 
 for $\phi$ sufficiently small.
Let $\Om_\phi$ be the family of $n-$forms which are
holomorphic in  the complex structure corresponding to $\phi$
and which are normalized according to the formula (\ref{eq:norm}).

We see that we get a family of (formal)
moduli spaces ${\M}(\phi)$ 
equipped
with 
coordinate systems $\mbox{Gr}\,W(\phi)\to 
{\M}(\phi)$, 
multiplications on $T_{{\M}(\phi)}$, symmetry vector fields,
and the pairings induced by  generalized period mappings from
the pairings defined by $\Om_\phi$.
We claim that this family of structures essentially does not depend on the
base point $\phi\in {\mathcal M}^{clasical}$.
The graded vector spaces $\mbox{Gr}\,W(\phi)$ are
identified by the Gauss-Manin connection.
Consequently one gets the identification of  the 
moduli spaces $\M(\phi)$ for  sufficiently small
 $\phi$.
Denote by $\tau^\phi$ the 
linear coordinate on $W_{\le n-1}/W_{\le n-\frac{3}{2}}$
 corresponding 
to $\phi$.

\begin{prp}\label{th:dtau=drho}
\begin{equation}
\frac{\p}{\p \phi} 
A_{ab}^{c}(\tau_{ W},{\phi})|_{\phi=0}=\frac{\p}
{\p \tau^\phi}A_{ab}^{c}(\tau_{W},0)
\end{equation}
\end{prp}
\proof
Let $\V(\phi), \V_{\hbar=1}(\phi)$ denote the differential graded Lie 
algebras associated with the base point $X_\phi$. Let 
$$
(\rho_\phi,f_\phi)\in\Om^{0,1}(X(0),\mathcal{T})\oplus \Om^{0,0}
(X(0),\mathcal{O})$$ denote the elements describing according to
the proposition \ref{pr:(ga,f)}
the deformations of the complex structure corresponding to 
$X_\phi$  and  the holomorphic volume element $\Om_\phi$
The proof follows from the standard quasi-isomorphisms
$\V(\phi)\sim (\V(0),d:=\db+[\rho_\phi,\cdot])$,
$\V_{\hbar=1}(\phi)\sim (\V_{\hbar=1}(0),
d:=\db+\dl+[\rho_\phi+f_\phi,\cdot])$.
\pendproof
It was shown in the proposition \ref{pr:PiMcl} that
\begin{equation}
\label{eq:Om(q,t)}
\Om(\phi)=\exp(\wg^W(\tau^\phi,\hbar=1))\vdash\Om(0)
\end{equation}
From this fact it can be deduced analogously that
\begin{equation}
\frac{\p}{\p\phi}(\Phi^W)^*(\p_a,\p_b)=0
\end{equation}

\subsection {Frobenius manifolds}

Let us recall  the   definition of formal  Frobenius (super) manifold
 as given 
in  \cite{km}, \cite{du}.  Let $\mathcal{H}$ be a finite-dimensional 
${\Z}_2$-graded affine space over\footnote{One can work over an
 arbitrary field of characteristic zero}
 $\C$. 
It 
is convenient to choose 
some set of coordinates $x_{\mathcal{H}}=\{x^a\}$ which defines the basis 
$\{\p_a:=\p/\p x^a\}$
of vector fields which are generators of affine transformations.
Let one of the coordinates is   distinguished and 
 denoted by $x_0$.
 Let $A^c_{ab}\in \C [[x_{\mathcal{H}}]]$ be a formal power series
  representing 3-tensor field,
$g_{ab}$ be a  nondegenerate  symmetric pairing on $\mathcal{H}$.

One can use the $A_{ab}^c$ in order to define a structure of 
$\C [[x_{\mathcal{H}}]]$-algebra
on $\mathcal{H}\otimes\C [[x_{\mathcal{H}}]]$, the (super)space of all
continious derivations of 
$\C[[x_{\mathcal{H}}]]$,
 with multiplication denoted by $\circ$:
\begin{equation} \p_a\circ\p_b:=\sum_c A_{ab}^c\p_c\end{equation}
One can use $g_{ab}$ to define the symmetric $\C[[x_{\mathcal{H}}]]$-pairing
on $\mathcal{H}\otimes\C[[x_{\mathcal{H}}]]$:
\begin{equation}\langle \p_a,\p_b\rangle:=g_{ab}\end{equation}

These data define the structure
of formal Frobenius manifold on $\mathcal{H}$ iff the following equations 
hold:
\begin{description}
\item[(1)] (Associativity/Commutativity/Identity)
\begin{equation}\forall \,a,b,c,d\,\,\,\,\,\,\sum_{e}A^e_{ab} A_{ec}^d 
=(-1)^{\bar a(\bar b+\bar c)}\sum_{e}A_{bc}^e 
A_{ea}^d \end{equation}
\begin{equation}
\forall \,a,b,c \,\,\,\,\,\, A^c_{ba}=(-1)^{\bar a\bar b}A_{ab}^c 
\label{eq:1b}\end{equation}
equivalently, $A^c_{ab}$ are  structure constants of a 
supercommutative
associative \break $\C[[x_{\mathcal{H}}]]$-algebra. It is 
required also that $\p_0$ is an identity element of this algebra.
 
\item [(2)](Invariance)  \begin{equation}\forall \,a,b,c 
\,\,\,\,\,\,\langle a\circ b,c\rangle=\langle a,b\circ 
c\rangle,\label{eq:2}\end{equation} 
equivalently,
the pairing $g_{ab}$ is invariant with respect to the 
multiplication~$\circ$. 

\item [(3)] (Potential) \begin{equation}\forall \,a,b,c,d\,\,\,\,\,\, \p_d 
A^c_{ab}=(-1)^{\bar a\bar d}\p_a 
A^c_{db}\,\,\,\, ,\end{equation}  which implies, assuming (\ref{eq:1b})
 and (\ref{eq:2}), that the series 
$A_{abc}$ are the third derivatives of a single power series
$\C[[x_{\mathcal{H}}]]$
\begin{equation}A_{abc}=\p_a\p_b\p_c\Phi\end{equation}
\item [(4)](Euler vector field)
There exists an even
 vector field $E$ acting as a conformal symmetry with respect to the
multiplication  and the metric  on ${{\mathcal T}}_{\mathcal{H}}$ in the 
following sense:
\begin{gather}\mbox{Lie}_E(g)=D g, \,\,\,D\in {\Bbb Q} 
\\
\mbox{Lie}_E(\circ)=\circ 
\end{gather} 
\end{description}
\subsection{Family of Frobenius manifold structures on $\M$}
\label{s:frob}
Recall that 
the set of rational Gromov-Witten invariants 
of a projective manifold $Y$ defines the Frobenius manifold
structure on $\oplus_{p,q}H^q(Y,\Om^p)[-q-p]$ 
considered as a (super)affine space.
Analogously, given a filtration $W$ complementary
to the Hodge filtration
the set of the invariants of generalized variations
of Hodge structures defines the Frobenius manifold structure
on  $\M$ identified with (super)affine space
$\oplus_{p,q}H^q(X,\Lambda^p\mathcal{T}_X)[-q-p]$ with the help
of the map 
$\mbox{Gr}\,\Phi^{W}$ (see \S \ref{s:invVHS}).

\begin{thm}\label{th:frob}
Assume that
we are given an isotropic  filtration $W$ 
on the graded vector space $H^*(X,\C)$
which is complementary
to the Hodge filtration. Then 
the affine structure on ${\M}$ corresponding to the 
series $\wg^W|_{\hbar=0}(\tau_W)$,
the
algebra structure on $T_{{\M}}$, the pairing
$(\Phi^W)^*(\,,\,)$ induced by the generalized period mapping
and the symmetry vector field $E$ 
define Frobenius
manifold structure on ${\M}$. 
\end{thm}
\proof
We have checked all the required properties above in this section.
\pendproof
\Remark
Frobenius manifold structure 
on $\M$ described in  \cite{bk}  is obtained if one takes as the
filtration $W$
the complex conjugate to the Hodge filtration $\overline{F}$.
\kendremark
\Remark
The expression (\ref{eq:Aabc}) suggests that there might exists a formula
of Chern-Simons type (see \cite{bk},\cite{bcov}) for the function
$\mathcal{F}^W$ whose third derivative equals $A_{abc}(\tau_W)$
\kendremark

\section{Mirror Symmetry in dimensions $n>3$}\label{ch:ms}
Let $(X,Y)$ be a pair of mirror dual Calabi-Yau manifolds. 
Hypothetically there exists singular point on the
compactification of the
moduli space of deformations of  $X$ such that 
the associated  limiting weight filtration $W$ is complementary 
to the Hodge filtration on $\oplus_k H^k(X,\C)[-k]$ and
its associated quotent $\mbox{Gr}\,W$ is naturally isomorphic 
to $\oplus_k H^k(Y,\C)[-k]$. Let us consider the
invariants of generalized variations of Hodge structure 
associated with $X$  which are written using the filtration $W$
at the singular point.
\vskip 3truemm

\noindent{\bf Mirror Symmetry conjecture in higher dimensions.}
{\it The  rational Gromov-Witten invariants of $Y$
coincides with the 
 invariants  (\ref{eq:inv}) of generalized
variations of Hodge structure associated with $X$.}
\vskip 2truemm

In this section we prove this conjecture in the case of projective complete
intersections.

\subsection{Projective hypersurfaces}
Let ${Y}^{n}\subset\Prj^{n+1}$ denotes a smooth
variety  defined by the equation 
$$P(x_1,\ldots,x_{n+2})=0$$
of degree $(n+2)$ in $(n+2)$ homogenious coordinates $(x_1:\ldots:x_{n+2})$.
The canonical class of ${Y}^{n}$ is trivial.
A nonzero section of the sheaf of holomorphic $n-$forms
 can be written as
\begin{equation}
\sum_{i=1}^{i=n+2}(-1)^i
\frac{x_idx_1\wedge\ldots\wedge\widehat{dx_{i}}\wedge\ldots\wedge 
dx_{n+2}}{dP}
\end{equation}
Let $X_{z}^{n}$ denotes the family of 
Calabi-Yau varieties obtained by the resolution of singularities
from the  varieties:
\begin{equation}
\bigl\{(x_1:\ldots:x_{n+2})|
x_1^{n+2}+\ldots+x_{n+2}^{n+2}=z\cdot x_1x_2\ldots 
x_{n+2}
\bigr\}/\bigl(\Z/(n+2)\Z\bigr)^{n}
\end{equation}
where
$\bigl(\Z/(n+2)\Z\bigr)^{n}$ is the group of diagonal matrices
acting on the zero locus of the above equation
\begin{equation}
\bigl\{\mbox{diag}(\zeta_1,\ldots,\zeta_{n+2})|\zeta^{n+2}_i=1,
\prod_{i=1}^{i=n+2}\zeta_i=1\bigr\}/\bigl\{
\mbox{diag}(\zeta,\ldots,\zeta)|\zeta^{n+2}=1\bigr\}
\end{equation}
The $1-$parameter family of varieties $X_z^n$ is mirror dual
to the universal family of Calabi-Yau hypersurfaces $Y^n\subset
\Prj^n$ (see\cite{bvs},\cite{gmp}).
The varieties $X^n$ and $Y^n$ are known to satisfy (see \cite{bb})
\begin{equation}\label{eq:hpq}
\mbox{dim} \,\,\,\,H^q(Y,\Om^p)=\mbox{dim} \,\,\,\,H^q(Y,\Om^{n-p})
\end{equation}

\subsection{Complete intersections}
More generally, let $Y_n(l_i)$ 
be a smooth variety defined as the intersection
of the hypersurfaces $\mathcal{L}_1,\ldots,\mathcal{L}_r\subset {\Bbb 
P}^{n+r}$ of degrees
$l_1,\ldots, l_r,\,\,l_1+\ldots+l_r=n+r+1$.
A construction of the mirror dual family of Calabi-Yau manifolds
$X^n_z(l_i)$ was proposed in \cite{bvs}.
The smooth varieties $X^n_z(l_i)$ are the Calabi-Yau 
 compactifications
of the complete intersection of $r$ hypersurfaces
\begin{equation}
\left\{ 
\begin{array}{lll}
u_1+&\dots&+u_{l_1}=1\\
u_{l_1+1}+&\dots&+u_{l_1+l_2}=1\\
&\dots& \\
u_{l_1+\dots+l_{r-1}}+&\dots&+u_{l_1+\dots+l_{r}}=1
\end{array}\right.
\end{equation}
inside the algebraic tori  
$\{u_1\cdot u_2\dots \cdot u_{n+r+1}=z|z\neq 0\}\subset {\Bbb A}^{n+r+1}$.
The pairs of varieties $Y_n(l_i),X_n(l_i)$ 
have Hodge numbers satisfying (\ref{eq:hpq})
according to \cite{bb}.

\subsection{Weight filtration}\label{s:W}
Let us denote by $\nu:\mathcal{X}(l_i)\to S$ the total
space of the family $X_z(l_i)$.
Let $W'_*$ denotes the standard limiting weight filtration 
on $\oplus_k H^k(X^n_z,\C)[-k]$ associated  with  the singular point 
$z=0$ in the family $X_z^n$. 
Recall that there exists a unique increasing filtration 
$W'_{\le 0}\subset W'_{\le 1}\subset\ldots\subset W'_{\le 2n}$
associated with
a nilpotent operator $N$ acting on a vector space $V$ with 
the properties
\begin{itemize}
\item $N:W_{\le l}'\subset W'_{\le l-2}$
\item $N^l:\mbox{Gr}_{n+l}W'\to\mbox{Gr}_{n-l}W'$ is
an isomorphism
\end{itemize}
The limiting weight filtration on $\oplus_k H^k(X_z,\C)$  
 is the  filtration defined by  the
nilpotent operator $N=\log(T)$ where $T$ is  the unipotent part of
the monodromy around $z=0$.  
Denote by  $W_{\le r}=W'_{\le 2r}$ the 
same filtration with all the indexes divided by two:
\begin{equation}
W_{\le 0}\subset W_{\le \frac{1}{2}}\subset\ldots \subset 
W_{\le n}=H^*(X,\Q)\end{equation}
It can be shown using the technic
from \cite{bb} and \cite{dl}
 that 
$W_*$ is complementary to the Hodge 
filtration $F^*(X_z)$ in the sense  of (\ref{eq:compl})
for sufficiently small $z$ and that the monodromy is in fact 
unipotent
(this is an extension of the  
condition of the (strong) maximal degeneracy of the point $z=0$,
see for example \cite{d},\cite{m}).

Let us introduce the following bigrading on
the  space $\mbox{Gr}\,W=\oplus_r W_{\le r}/W_{\le r-\frac{1}{2}}$:
\begin{multline}
\mbox{Gr}\,W=\oplus_{i,j\in 0\ldots n}
(\mbox{Gr}\,W)^{ij},\\
(\mbox{Gr}\,W)^{ij}=\bigl(W_{\le n-(i+j)/2}\cap
H^{j-i+n}(X,\C)\bigr)/\bigl(W_{\le n-(i+j+1)/2}\cup
H^{j-i+n}(X,\C)\bigr)
\label{eq:bigr}
\end{multline}
Since  $W$ is complementary to the Hodge filtration $F$
one has  $\mbox{Gr}\,W\simeq\mbox{Gr}\, F$ and
\begin{equation}
(\mbox{Gr}\,W)^{ij}\simeq H^{j}(X,\Om^{n-i}_X)
\end{equation}

A choice of holomorphic $n-$form $\Om_X$ defines a pairing (\ref{eq:pair}) 
on the space $\mbox{Gr}\,F^H=\oplus_{i,j}H^{j}(X,\Lambda^iT_X)[-i-j]$.
It induces a natural pairing  on the
the space $\mbox{Gr} \,F=\oplus_{i,j} H^{j}(X,\Om^{n-i}_X)[-i-j]$ 
via the identification $\mbox{Gr}\,F^H\stackrel{\vdash\Om_X}{=}\mbox{Gr}\,F$
where the same choice of the holomorphic
$n-$form $\Om_X$ is used. 
When considered as a pairing on $\mbox{Gr}\,W\simeq\mbox{Gr}\,F$ 
this pairing does not depend either on 
the deformations of complex structure on $X$ nor on the choice
of $\Om$.
This pairing is symmetric with respect to
the    total grading (\ref{eq:bigr}) on $\mbox{Gr}\,W$.
It differs from the pairing
induced from the Poincare pairing by certain signs depending
on the given bigraded component (see eq.~(\ref{eq:dg=0})).

It is known that the operator $N$ acting 
on $\mbox{Gr}\,W$ has 
properties similar  to that of the
operator of multiplication by the K\"ahler 2-form 
acting on the cohomologies of a projective manifold.
In  particular, the operator $N$ defines
  the ``Lefschetz decomposition''
on $\mbox{Gr}\,W$.
The restriction of the pairing to the corresponding primitive subspaces 
satisfies the analogs of the Hodge-Riemann bilinear identities (see 
\cite{d} and references therein).
 
Recall that the only nonzero Hodge numbers of $Y$ are 
$h^{p,q}$ for $p+q=\mbox{dim}_{\C}Y$, and $h^{p,p}=1\,\,
\mbox{for}\,\,p\neq \mbox{dim}_{\C}Y/2$. 
Let us denote by $L$ the operator of the
multiplication by the K\"ahler 2-form acting on $H^*(Y,\C)$.
The Lefschetz decomposition on $H^*(Y,\C)$ has the following
structure
\begin{gather}
H^*(Y,\C)=<\Delta_0,\ldots,\Delta_n>\oplus(\oplus_{k=1\ldots\mbox{dim}\,
P^n(Y,\C)}<\Theta_k>)\label{eq:basis}\\
\Delta_0=1,\,\,\,L\Delta_i=\Delta_{i+1},\, 
\,\,L\Delta_n=0,\,\,\,\Delta_i\in 
H^{2i}(Y,\C),\,\,\,\,n=\mbox{dim}_{\C}Y\notag\\
L\Theta_k=0,\,\,\,\,\Theta_k\in 
H^{n}(Y,\C)\notag\end{gather}
It follows from the property (\ref{eq:compl}) that the
"Lefschetz decomposition" 
induced by the operator $N$
on $\mbox{Gr}\,W$ has similar 
structure.
\begin {prp}\label{pr:f}
There exists an  isomorphism of the $\Z-$graded
\footnote{When we speak about the $\Z-$grading 
on $\mbox{Gr}\,W$ we mean the total sum of gradings
(\ref{eq:bigr})} vector spaces
\begin{equation}
f:H^*(Y,\C)\to \mbox{Gr}\,W
\end{equation}
preserving the  pairing and such that
$f^{-1}Nf=L$
\end{prp}
\proof It follows from the coincidence of the 
dimensions of the components  of the Lefschetz decompositions.
\pendproof
Let us denote 
\begin{equation}
\widetilde\Delta_i=f(\Delta_i),\,\,\,\,\,
\widetilde\Theta_k=f(\Theta_k)\label{eq:Xbasis}
\end{equation}

We will use the linear map $f$ in \S \ref{s:supermain}
to compare the Gromov-Witten invariants
of $Y^n(l_i)$ with the invariants of the generalized
variations of Hodge structures attached to $X^n(l_i)$ 

\subsection{q-expansion}
Let 
\begin{gather}
\G_i\in H_n(X_z,\C)\cap W_{\le -i}^{\bot},\,\,\,\,
N\G_i=
\G_{i-1} \notag\\
\Xi_k\in W_{\le-n/2}^{\bot},\,\,\, N\Xi_k=0
\label{eq:Gi}
\end{gather}
the elements which form a frame in $H_*(X_z,\C)$ 
which projects to the   frame in $\mbox{Hom}(\mbox{Gr}\,W,\C)$ dual to
$\{\widetilde\Delta_i,\widetilde\Theta_k\}$.
It follows that the monodromy transforms
 $\G_0, \G_1$ into
\begin{equation}
T(\G_0)=\G_0,\,\,\,T(\G_1)=\G_1+\G_0\label{eq:TG1}
\end{equation}
Let $\Om(z)$ denotes the holomorphic
$n-$form on $X_z$ normalized so that
\begin{equation}
\int_{\G_0}\Om(z)=1\label{eq:omg0}
\end{equation}
Recall that 
we have associated the formal power series
$A_{abc}(\tau,z)$ to the choice of
the base point $z$ of the classical moduli space of complex structures,
of the isotropic filtration $W$ complementary to the Hodge filtration,
and of the normalization of the holomorphic $n-$form $\Om(z)$.

According to
the formula (\ref{eq:coor})
affine coordinate induced on the base of the family $X_z^n$
via the map 
$\mbox{Gr}\,\Phi^W|_{\M^{classical}}$
can be identified with the period 
\begin{equation}
t^{base}(z)=\int_{\G_1}\Om(z) \label{eq:omg1}
\end{equation}
It follows from the formula (\ref{eq:TG1}) that
the family of formal power series
$A_{abc}(\tau;t^{base}(z))$ depends on $t^{base}(z)$ only
through $q=\exp((2\pi i)t^{base}(z))$. Denote by $\tau^{\G_1}$ the 
coordinate on the formal moduli space $\M(q)$ corresponding to the
element $[\G_1]\in (W_{\le n-1}/W_{\le n-\frac{3}{2}})^{dual}$. 
\begin{prp}\label{pr:a(qet)}
\begin{equation}
A_{abc}(\tau,\tau^{\G_1};q)=
A_{abc}(\tau, 0;\exp(2\pi i\,\tau^{\G_1})q)
\end{equation}
considered as  power series in $(\tau,q)$
\end{prp}
\proof
It follows immediately
from the prop.~\ref{th:dtau=drho}.
\pendproof
One can assume without loss of generality that
the element $\G_1$ defined up to an addition of 
an element proportional to $\G_0$ is chosen so that
the coordinates $q$ on the base of the family $X(l_i)$ coincides
with the similar  coordinate from \cite{bvs},\cite{g1}.

\subsection{Invariants of the generalized VHS and Gromov-Witten invariants}
\label{s:supermain}
Let $C_{abc}(\tau',q),\,\,\tau'\in H^*(Y^n,\C)$ 
is the generating
function encoding the whole set of rational Gromov-Witten invariants
(for the definition of $C_{abc}(\tau',q)$ see \cite{km,bm}) 
of the smooth projective Calabi-Yau variety $Y^n(l_i)$ which is
the intersection of the hypersurfaces of degrees $l_1,\ldots,l_r$. 
We would like to identify $C_{abc}(\tau',q)$ with
the generating function  $A_{abc}(\tau,q)$
of the invariants introduced in section \ref{ch:inv} which are
attached to the dual family $X^n(l_i)$.

\subsubsection*{Picard-Fuchs equations}
Here we explain how to write down the equations satisfied by the
 periods of the holomorphic $n-$form $\Om(q)$:
\begin{equation}
\int_{\G_i}\Om(q)\label{eq:prds}
\end{equation}
where $N\G_i=\G_{i-1}$ are the locally constant elements defined in 
(\ref{eq:Gi}).

\newcommand{\wde}{\widetilde\Delta}

Let us denote  
\begin{gather}
\widetilde\Delta_i^q\in W_{\le n-i}\cap F^{\ge n-i}\cap H^n(X_q,\C),
\,\,\,\,\,<\wde_i^q,\ldots,\wde_{n}^q>=\im\,N^i,\notag\\
\widetilde\Theta_k^q\in W_{\le n/2}\cap F^{\ge n/2}
\end{gather}
the  
sections  of $R^*\nu_*(\C_{\mathcal{X}})$ whose images in $\mbox{Gr}\,W$ 
coincide with $\widetilde\Delta_i,\widetilde\Theta_k$.
\begin{prp}
The covariant derivatives of $\wde_i^q$ with respect to the Gauss-Manin
connection have the following form:
\begin{gather}
\mathcal{D}_q \wde_{i-1}^q=\frac{1}{q} a_i(q)\wde_{i}^q\,\,
\fr \,\,i\in 1\ldots n-1\label{eq:DqDelta}
\\ a_i(q)=1+o(q)\label{eq:ai(q)}
\end{gather}
\end{prp}
\proof
The operator $N$ is locally constant.
It follows that
the subspace generated by $\wde_i^q,\ldots,\wde_{n}^q$ 
is preserved by the Gauss-Manin connection.
The equation  (\ref{eq:DqDelta}) follows now from the
Griffiths transversality condition $D_q F^{\ge s}\subset F^{\ge s-1}$.
The form (\ref{eq:ai(q)}) of the coefficients $a_i(q)$ follows from the
nilpotent orbit theorem (see \cite{d} and references therein).
\pendproof

One can  check that $\Om(q)=\wde_0^q$. The conditions
(\ref{eq:omg0}) and (\ref{eq:omg1}) imply also that
\footnote{In the sequel we will denote
the coordinate $t^{base}$ simply by $t$ when it does not
seem to lead to a confusion} $\mathcal{D}_t\Om(q)=\wde_1^q$.
\begin{cor}
The periods  of the holomorphic $n-$form $\Om(q)$ satisfy
\begin{equation}
\p_t(\frac{1}{a_n(e^{2\pi it})}\p_t(\ldots\p_t(\frac{1}{a_1(e^{2\pi it})}
\p_t \,\,\int_{\G_i}\Om(e^{2\pi it})\,\,\,)\ldots)=0\label{eq:a(q)}
\end{equation}
\end{cor}
\pendproof
Let  
$$
\widetilde\Delta_i^H,\widetilde\Theta_k^H\in 
\oplus_{r,s}H^{r}(X_q,\Lambda^sT_{X_q})[s-r]
$$
correspond to $\widetilde\Delta_i,\widetilde\Theta_k\in \mbox{Gr}\,W$ 
via the
composition of isomorphisms 
$$
\mbox{Gr}\,W\simeq\mbox{Gr}\,F
\stackrel{\vdash\Om(q)}{\simeq}\oplus_{r,s}H^{r}(X_q,\Lambda^sT_{X_q})[s-r]
$$
In particular, $\wde_0^H=1$ and
 $\wde_1^H=[\frac{\p}{\p t}]$.
\begin{prp}
The following
identities hold 
\begin{equation}
\wde_1^H\cdot\wde_{i-1}^H=a_i(q)\wde_i^H\fr \,\,i\in 1\ldots n-1
\end{equation}
in the algebra 
$\oplus_{r,s} 
H^{r}(X_q,\Lambda^sT_{X_q})[s-r]$ 
\end{prp}
\proof
The covariant derivative of the flag $F_t$ is 
a linear map $$\mbox{Gr}\, F^{\ge *}_t\to 
\mbox{Gr}\, F_t^{\ge *-1}$$ 
It coincides 
as an element of
$$\oplus_{r,s}\mbox{Hom}\,(H^r(X_t,\Om^s_X),H^{r+1}(X,\Om^{s-1}_X))$$
with the cup multiplication by $\p/\p t$ 
\pendproof

\subsubsection*{Small quantum cohomology differential equations}
Here we decribe the analogous differential equation  involving
the structure constants of the small quantum multiplication
on $H^*(Y,\C)$. 

\begin{prp}
The following identities hold 
in the small quantum cohomology algebra of $Y(l_i)$ 
\begin{gather}
\Delta_1\cdot\Delta_{i-1}=c_i(q)\Delta_i\fr \,\,i\in 1\ldots n-1\\
c_i(q)=1+o(q)\label{eq:c(0)}
\end{gather}
\end{prp}
\proof
It follows almost immediately from the definitions.
The only thing to check which is not immediately obvious is the fact that
$\Delta_1\cdot\Delta_{\frac{n}{2} -1}\in\,\,<\Delta_{\frac{n}{2}}>$
in the case $\mbox{dim} \,Y$ is even.
This follows  from the fact that the calculation of the
small quantum multiplication $\Delta_1\cdot\Delta_{\frac{n}{2} -1}$
can be reduced to the calculation of certain intersection numbers
on the space of stable maps to the projective space 
$\C\mathbb{P}^{n+r+1}\supset Y(l_i)$.
\pendproof

Consider the differential equation of the $n+1$-st order
\begin{equation}
\p_t(\frac{1}{c_n(e^{2\pi it})}\p_t(\ldots\p_t(\frac{1}{c_1(e^{2\pi it})}
\p_t \,\,\psi(t)\,\,\,)\ldots)=0\label{eq:c(q)}
\end{equation}

\begin{prp}\label{pr:a(q)=c(q)}
For the pair of mirror families 
$X(l_i),Y(l_i)$ one has 
\begin{equation}
a^X_i(q)=c^Y_i(q)\label{eq:a(q)=c(q)}
\end{equation}
\end{prp}
\proof
The theorem 11.8 from \cite{g1} implies that
the two differential equations (\ref{eq:a(q)})
and (\ref{eq:c(q)}) have the same space of solutions.
This identifies $a^X_i$ and $c^Y_i$ up to
a multiplication by a constant.
The latter ambiguity is fixed by
the conditions (\ref{eq:ai(q)}) and (\ref{eq:c(0)}).
\pendproof

\subsubsection*{The two generating functions}
Here we prove the  theorem establishing the
coincidence of the Gromov-Witten invariants  
of the projective 
complete intersection Calabi-Yau manifolds
and the invariants introduced in the section \ref{ch:inv}
associated  with  their
mirror duals. Recall that we use the map 
$f:H^*(Y,\C)\to \mbox{Gr}\,W$ defined in the proposition
\ref{pr:f} 
to compare the generating power series $A(\tau,q)$ and $C(\tau,q)$.
 \begin{thm} 
\begin{equation}
C_{abc}(\tau,q)=A_{abc}(\tau,q)
\label{eq:a=c}\end{equation}
in other words the  rational Gromov-Witten invariants
of $Y^n(l_i)$ coincide with the  invariants of generalized
variations of Hodge structures associated with $X^n(l_i)$
\end{thm}
\proof
The idea of the proof is to use the constraints
on the series $A_{abc}(\tau,q)$, $C_{abc}(\tau,q)$ arising
from the equations defining the Frobenius manifold
structure and the proposition \ref{pr:a(q)=c(q)}. 
Let us consider the Taylor expansions
\begin{gather}
C_{abc}(\tau,q)=\sum_{m\geq 0, d_1\ldots d_k}
\frac{1}{k!}\varepsilon(a,b,c,d_1\ldots d_k)
<e_a,e_b,e_c,e_{d_1},\ldots,e_{d_k}>^C_m q^m \tau^{d_1}\ldots\tau^{d_k}
\\
A_{abc}(\tau,q)=\sum_{m\geq 0, d_1\ldots d_k}
\frac{1}{k!}\varepsilon(a,b,c,d_1\ldots d_k)
<e_a,e_b,e_c,e_{d_1},\ldots,e_{d_k}>^A_m q^m 
\tau^{d_1}\ldots\tau^{d_k}
\notag
\end{gather}
where $\{e_i\}$ denotes the basis  (\ref{eq:basis}), (\ref{eq:Xbasis})
in
$H^*(Y,\C)\stackrel{f}{\simeq} \mbox{Gr}\, W$ and
$\varepsilon$ is the standard sign depending on the
parity of the elements $e_{d_i}$.
Recall 
 that the conformal symmetry vector field
acts on both power series
(see eq.~(\ref{eq:Ecirc}) for the case of $A_{abc}$)
\begin{equation}
E(\tau)=\sum_a 
\frac{1}{2}(\mbox{deg}\,\tau^a+2)\tau^a\frac{\p}{\p\tau^a}
\end{equation}
where $\mbox{deg}\,\tau^a=-k$ for $\tau^a\in 
(W_{\leq (n-k)/2}/ W_{\leq (n-k-1)/2})^{dual}
\stackrel{f^*}{\simeq}(H^{k}(Y,\C))^{dual}$.
It follows that
\begin{equation}
<e_{a_1},\ldots,e_{a_k}>^{A(C)}_m\neq 0\Rightarrow
\sum_i \mbox{deg}\,e_{a_i}=2\mbox{dim}_\C X+ 2(k-3) \label{eq:grading}
\end{equation}
where $\mbox{deg}\,e_a=k$ for $e_a\in 
H^{k}(Y,\C)\stackrel{f}{\simeq}
W_{\leq (n-k)/2}/ W_{\leq (n-k-1)/2}$.
It follows from the proposition \ref{pr:1}
that $A_{ab0}(\tau,q)=\eta_{ab}$ where $\eta_{ab}$ is
the $2-$tensor of the metrics.
The proposition \ref{pr:a(q)=c(q)}  
together with the grading conditions defined by the
symmetry vector field $E(\tau)$
  imply that
$C_{abc}(0,q)=A_{abc}(0,q)$
whenever one of the indexes $a,b,c$ corresponds
to $H^2(Y,\C)\simeq (\mbox{Gr}\,W)^2$.
The proposition \ref{pr:a(qet)} implies that the
following
analog of the "Divisor axiom" holds for the
series $A_{abc}(\tau,q)$
\begin{multline}
\forall\,\,\, a,b,c,d_1\ldots d_k; m\\
<e_a,e_b,e_c,e_{d_1},\ldots,e_{d_k},\wde_1>^A_m=
<e_a,e_b,e_c,e_{d_1},\ldots,e_{d_k}>^A_m
\end{multline}
where $\wde_1$ is the basis element corresponding to ${\p}/
{\p t}^{base}$.
Let us consider the associativity equaion
\begin{equation}
\forall \,\,a,b,c,d,\,\,
\sum_{f,g} A_{abf}\eta^{fg}A_{gcd}=(-1)^{\bar a(\bar b+\bar c}
\sum_{f,g}A_{bcf}\eta^{fg}A_{gad}\label{eq:assoc}
\end{equation}
and similarly for $C_{abc}$.
Notice it follows from the grading condition (\ref{eq:grading})
that any  
nonzero expression $<e_{i_1},\dots,e_{i_k}>^{C}_m$ 
contains no more then two elements from the ``nonalgebraic'' subspace 
generated by $\Theta_i$. Analogously the same is 
true for the Taylor 
coefficients of the series $A_{abc}(\tau,q)$ and the elements from the 
subspace generated by
$\widetilde \Theta_i$. 
Using the equation (\ref{eq:assoc}) as in the proof of
 the theorem 3.1 from \cite{km} 
  all the Taylor coefficients
 of the series $A_{abc}$ and $C_{abc}$ can now be identified inductively.
\pendproof

Let $\{\mathcal{G}_i\}$ denotes a locally constant frame
in $H_*(X_q,\C)$.
\begin{cor}
The Gromov-Witten invariants of $Y(l_i)$
are expressed in terms of the generalized periods $\int_{\mathcal{G}_i}\Pi^W
(\tau,q)$
(see formula (\ref{eq:exp}))
associated
with the dual family $X(l_i)$ 
\begin{equation}
C_{ab}^{c}(\tau,q)=
\sum_i((\p\Pi)^{-1})^c_i\p_{a}\p_{b}\Pi^i
\end{equation}
\end{cor}





\noindent{Address:\\ S.Barannikov}

\noindent{IHES}

\noindent{35, route de Chartres}

\noindent{91440 France}

\vskip 3truemm
\noindent{e-mail: barannik@ihes.fr}

\end{document}